\documentclass[12pt]{article}

\usepackage[latin1]{inputenc}
\usepackage{bbm}
\usepackage{stmaryrd}
\usepackage{latexsym}
\usepackage{amsfonts}
\usepackage{amsmath,amssymb,amscd}
\usepackage{yfonts}
\usepackage{dsfont}
\usepackage{mathabx}
\usepackage[dvips]{graphicx}
\usepackage{epsfig}
\usepackage{psfrag}
\usepackage[font={small}]{caption}

\DeclareFontFamily{U}{txsyc}{}
\DeclareFontShape{U}{txsyc}{m}{n}{
   <-> txsyc%
}{}
\DeclareFontShape{U}{txsyc}{bx}{n}{
   <-> txbsyc%
}{}
\DeclareFontShape{U}{txsyc}{l}{n}{<->ssub * txsyc/m/n}{}
\DeclareFontShape{U}{txsyc}{b}{n}{<->ssub * txsyc/bx/n}{}
\DeclareSymbolFont{symbolsC}{U}{txsyc}{m}{n}
\SetSymbolFont{symbolsC}{bold}{U}{txsyc}{bx}{n}
\DeclareFontSubstitution{U}{txsyc}{m}{n}
\DeclareMathSymbol{\df}{\mathrel}{symbolsC}{"42}
\DeclareMathSymbol{\fd}{\mathrel}{symbolsC}{"43}
\DeclareMathSymbol{\lJoin}{\mathrel}{symbolsC}{"58}
\DeclareMathSymbol{\rJoin}{\mathrel}{symbolsC}{"59}

\newcommand{\cA}{{\cal A}}

\newcommand{\cC}{{\cal C}}

\newcommand{\cH}{{\cal H}}

\newcommand{\cK}{{\cal K}}
\newcommand{\cL}{{\cal L}}

\newcommand{\cP}{{\cal P}}

\newcommand{\cR}{{\cal R}}
\newcommand{\cS}{{\cal S}}
\newcommand{\cT}{{\cal T}}

\newcommand{\CC}{\mathbb{C}}

\newcommand{\EE}{\mathbb{E}}

\newcommand{\LL}{\mathbb{L}}
\newcommand{\NN}{\mathbb{N}}
\newcommand{\PP}{\mathbb{P}}

\newcommand{\RR}{\mathbb{R}}
\renewcommand{\SS}{\mathbb{S}}

\newcommand{\ZZ}{\mathbb{Z}}

\newcommand{\iy}{\infty}
\newcommand{\lt}{\left}
\newcommand{\me}{\medskip}

\newcommand{\pa}{\partial}
\newcommand{\ri}{\rightarrow}
\newcommand{\rt}{\right}
\newcommand{\sm}{\smallskip}
\newcommand{\tr}{\triangle}

\newcommand{\wi}{\widetilde}
\newcommand{\wit}{\widehat}

\newcommand{\ex}{\exists\ }
\newcommand{\fo}{\forall\ }

\newcommand{\lve}{\lt\vert}

\newcommand{\rve}{\rt\vert}

\newcommand{\st}{\,:\,}

\newcommand{\un}{\mathds{1}}

\newcommand{\bq}{\begin{eqnarray*}}
\newcommand{\bqn}[1]{\begin{eqnarray}\label{#1}}
\newcommand{\eq}{\end{eqnarray*}}
\newcommand{\eqn}{\end{eqnarray}}
\newcommand{\wwtbp}{\par\hfill $\blacksquare$\par\me\noindent}
\newcommand{\thistitlepagestyle}{}
\newcommand{\lin}{\llbracket}
\newcommand{\rin}{\rrbracket}

\newcommand{\ttsim}{\raise.17ex\hbox{$\scriptstyle\mathtt{\sim}$}}

\newtheorem{pro}{Proposition} 
\newtheorem{cor}[pro]{Corollary}
\newtheorem{lem}[pro]{Lemma}
\newtheorem{theo}[pro]{Theorem}
\renewcommand{\thepro}{\arabic{pro}}

\newenvironment{rem}
{\par\me\refstepcounter{pro}\noindent{\bf Remark \thepro\ }}
{\par\hfill $\square$\par\sm\noindent}
\newenvironment{rems}
{\par\me\refstepcounter{pro}\noindent{\bf Remarks \thepro\ }\par\sm}
{\par\hfill $\square$\par\sm\noindent}
\newcommand{\proof}{\par\me\noindent\textbf{Proof}\par\sm\noindent}
\newcommand{\prooff}[1]{\par\me\noindent\textbf{#1}\par\sm\noindent}

\newcommand{\comment}[1]{}

\title{On the Markov commutator}
\author{Laurent Miclo
}
\date{\box1
}

\begin{document}

\setbox1=\vbox{
\large
\begin{center}
Institut de Mathématiques de Toulouse, UMR 5219\\
Université de Toulouse and CNRS, France\\
\end{center}
} 
\setbox3=\vbox{
\hbox{miclo@math.univ-toulouse.fr\\}
\vskip1mm
\hbox{Institut de Mathématiques de Toulouse\\}
\hbox{Université Paul Sabatier\\}
\hbox{118, route de Narbonne\\} 
\hbox{31062 Toulouse Cedex 9, France\\}
}
\setbox5=\vbox{
\box3
}

\maketitle
\thistitlepagestyle
\abstract{
The Markov commutator associated to a finite Markov kernel $P$ is
the convex semigroup consisting of all Markov kernels commuting with $P$.
Its interest comes from its relation with the hypergroup property and with
the notion of Markovian duality by intertwining.
In particular, it is shown that the discrete analogue of the Achour-Trimèche's theorem, asserting the preservation of non-negativity
by the wave equations associated to certain Metropolis birth and death transition kernels, cannot be extended
to all convex potentials. But it remains true for symmetric and monotone potentials which are sufficiently convex.
}
\vfill\null
{\small
\textbf{Keywords: }
finite Markov kernels, Markov commutator, symmetry group of a Markov kernel, hypergroup property, duality by intertwining, Achour-Trimèche theorem, birth and death chains, Metropolis algorithms, one-dimensional discrete wave equations.
\par
\vskip.3cm
\textbf{MSC2010:} primary: 
 60J10, secondary: 
 15A27, 20N20, 52C99, 39A12.
}\par

\newpage

\section{Introduction}

The primary motivation for this paper is to disprove, at least in a finite context, a conjecture due to Dominique Bakry,
about an extension of Achour-Trimèche's theorem \cite{MR552062} (see also Bakry and Huet \cite{MR2483738}).
It also provides the opportunity to begin a systematic study of the commutator convex semi-group associated to a Markov kernel.\par\me
Here we will only  be concerned with  state spaces $V$ which are finite and endowed with a Markov kernel $P$, namely
a matrix $(P(x,y))_{x,y\in V}$ whose entries are non-negative and whose row sums are equal to 1.
Two classical assumptions on $P$ are:\\
\textbf{Irreducibility}: all the coefficients of 
$\sum_{n\in\lin \lve V\rve\rin} P^n$
are positive ($\lve V\rve$ is the cardinality of $V$ and we denote for any $k\leq l\in\ZZ$, $\lin k,l\rin\df \{k, k+1, ..., l-1, l\}$,
and $\lin k\rin\df \lin 1,k\rin$ for $k\in\NN$).\\
\textbf{Reversibility}: there exists a probability measure $\mu$ positive on $V$, such that
\bqn{rev}
\fo x,y\in V,\qquad \mu(x)P(x,y)&=&\mu(y)P(y,x)\eqn
\par
Under the reversibility assumption, there exist orthonormal bases of $\LL^2(\mu)$
consisting of eigenvectors $\varphi_1,\ \varphi_2, ...,\ \varphi_{\lve V\rve}$ of $P$, 
associated to the eigenvalues $1= \theta_1\geq \theta_2\geq \cdots\geq \theta_{\lve V\rve}\geq -1$.
Without loss of generality, we will always choose $\varphi_1=\un$.
We say that $P$ satisfies the \textbf{hypergroup property} with respect to a point $x_0\in V$, if the previous basis
can be chosen such that $\varphi_k(x_0)\not=0$ for all $k\in\lin\vert V\vert\rin$, and
\bqn{Bakry}
\fo x,y,z\in V,\qquad \sum_{k\in\lin \vert V\vert\rin} \frac{\varphi_k(x)\varphi_k(y)\varphi_k(z)}{\varphi_k(x_0)}&\geq &0\eqn
\par
These notions can be immediately extended to Markov generators $L$ on $V$, namely matrices whose off-diagonal entries are non-negative
and whose row sums vanish (for instance by considering the generated semi-group $(P_t)_{t\geq 0}\df (\exp(tL))_{t\geq 0}$
and by asking that the above conditions are satisfied by $P_t$, for some $t>0$, it does not depend on the choice of $t>0$).
Extensions to more general Markov processes are also possible, but they may require some care.
E.g.\ in \cite{MR2483738}, Bakry and Huet consider one-dimensional diffusion generators of the form $L_U\df \pa^2-U'\pa$ on $[-1,1]$,
with Neumann conditions on the boundary and where $U\st [-1,1]\ri \RR$ is a smooth potential.
They prove Achour-Trimèche's theorem \cite{MR552062}, asserting that if $U$ is convex and either
monotonous or symmetric with respect to 0, then $L_U$ satisfies the hypergroup property.
In a personal communication, Dominique Bakry was wondering if this result would remain true if the assumption
``monotonous or symmetric with respect to 0" was removed.
Our main objective is to show that this is wrong, at least in the finite setting.
\par
More precisely, let $N\in\NN\setminus\{1\}$ be given and denote by $\cC$ the set of functions $U\st \lin 0, N\rin\ri \RR$
which are convex (i.e.\ whose natural piecewise affine extension to $[0,N]$ is convex).
For $U\in \cC$, let $\mu_U$ be the probability on $\lin 0,N\rin$ given by
\bqn{muU}
\fo x\in \lin 0, N\rin,\qquad
\mu_U(x)&\df& Z_U^{-1}\exp(-U(x))\eqn
where $Z_U$ is the renormalizing constant.
For any $U\in\cC$, assume we are given an irreducible birth and death Markov transition $P_U$ on $\lin 0,N\rin$ whose invariant probability 
is $\mu_U$. Recall that a \textbf{birth and death} Markov transition $P$ on $\lin 0,N\rin$ is a Markov kernel such that
\bq
\fo x,y\in \lin 0,N\rin,\qquad P(x,y)>0&\Rightarrow& \vert x-y\vert \leq 1\eq
An invariant measure of such a kernel necessarily satisfies \eqref{rev}, so that an irreducible birth and death Markov matrix is reversible.
\par
Endowing $\cC$ and the set of Markov kernels from the topology inherited respectively from $\RR^{\lin 0, N\rin}$ and $\RR^{\lin 0, N\rin^2}$,
we say that the above mapping $\cC\ni U\mapsto P_U$ is a (birth and death) \textbf{generalized Metropolis procedure} if it is continuous.
A classical Metropolis procedure corresponds for instance to the Markov kernel $M_U$ defined by
\bqn{MU}
\fo x\not=y\in\lin 0,N\rin,\qquad M_U(x,y)&\df& \frac{M_0(x,y)}{\Sigma_U}\exp\lt(\frac{U(x)-U(y)}{2}\rt)\eqn
where the exploration Markov kernel $M_0$ is given by
\bqn{M0}
\fo x\not=y\in\lin 0,N\rin,\qquad M_0(x,y)&\df&\lt\{
\begin{array}{ll}
1/2&\hbox{, if $\lve x-y\rve=1$}\\
0&\hbox{, otherwise}
\end{array}\rt.
\eqn
and where 
\bqn{SigmaU}
\Sigma_U&\df&\max_{x\in \lin 0,N\rin}\sum_{y\in\lin 0,N\rin\setminus\{x\}}M_0(x,y)\exp\lt(\frac{U(x)-U(y)}{2}\rt) \eqn
As usual, the diagonal entries of the matrices $M_U$ and $M_0$ are imposed by the condition that the row sums are equal to 1.
\par
Our main result is:
\begin{theo}\label{theo1}
It does not exist a generalized Metropolis procedure $\cC\ni U\mapsto P_U$
such that $P_U$ satisfies the hypergroup property for all $U\in\cC$.
\end{theo}
\par
In  \cite{miclo:hal-01117051}, we checked numerically (by appropriate random choices of $U$  in $\cC$) that a variant of the classical Metropolis procedure (described as $\cC\ni U\mapsto \wideparen{M}_U$ with the notation introduced in \eqref{wMU} below)
does not satisfy the hypergroup property.
\par\sm
The proof of Theorem \ref{theo1} is based on properties of the \textbf{commutator convex semi-group} $\cK(P)$ associated to a Markov kernel $P$ on $V$:
it is the set of Markov kernels 
$K$ on $V$ commuting with $P$: $KP=PK$. It is immediate to see that it is convex and that it is a semi-group: if $K$ and $K'$ belong to $\cK(P)$,
the same is true for their product $KK'$.
It was introduced in \cite{miclo:hal-01117051}, because it gives a simple Markovian characterization of the hypergroup property for certain kernels.
More precisely, let us introduce the following objects:
\bq
\fo x\in V,\qquad \cK(P,x)&\df& \{K(x,\cdot)\st K\in\cK(P)\}\ \subset \ \cP(V)\eq
where $\cP(V)$ is the convex set of probability measures on $V$, and
\bq
\cH(P)&=&\{x\in V\st \cK(P,x)=\cP(V)\}\eq
Furthermore, say that a Markov kernel is \textbf{uniplicit} if it is reversible and if all its eigenvalues are of multiplicity 1
(in particular the eigenvalue 1 is of multiplicity 1, so that  uniplicity implies irreducibility).
The interest of these notions is:
\begin{lem}\label{cHgroup}
An uniplicit Markov kernel $P$ on $V$ satisfies the hypergroup property with respect to $x_0\in V$  if and only if $x_0\in\cH(P)$.
\end{lem}
Let us give  succinctly some underlying arguments, since this is the only place in the paper where Definition \eqref{Bakry} will play a role.
\proof
The reverse implication was observed in \cite{miclo:hal-01117051} and the direct implication is a consequence of the considerations of 
Bakry and Huet \cite{MR2483738}, the uniplicit assumption is not even needed, as the following reminder show. Let $P$ be a reversible Markov kernel $P$ on $V$ 
with an associated orthonormal basis of  eigenvectors
 $\varphi_1,\ \varphi_2, ...,\ \varphi_{\lve V\rve}$ as above.
Assume that $P$ satisfies the hypergroup property with respect to $x_0\in V$.
Let $x\in V$ be given and consider the kernel $K_x$ given by
\bq
\fo y,z\in V,\qquad K_x(y,z)&\df& \sum_{k\in\lin \vert V\vert\rin} \frac{\varphi_k(x)\varphi_k(y)\varphi_k(z)}{\varphi_k(x_0)}\, \mu(z)\eq
By assumption it is non-negative and for any fixed $y\in V$,  we have by orthonormality,
\bq
\sum_{z\in V}K_x(y,z)&=&\sum_{z\in V}K_x(y,z)\varphi_1(z)\\
&=& \sum_{k\in\lin \vert V\vert\rin} \frac{\varphi_k(x)\varphi_k(y)}{\varphi_k(x_0)}\sum_{z\in V}\varphi_k(z)\varphi_1(z)\, \mu(z)\\
 &=& \frac{\varphi_1(x)\varphi_1(y)}{\varphi_1(x_0)}\\
 &=&1\eq
 Thus $K_x$ is a Markov kernel. A similar computation shows that for any $k\in\lin 2,\lve V\rve\rin$, $\varphi_k$ is also an eigenfunction of $K_x$ associated
 to the eigenvalue $\varphi_k(x)/\varphi_k(x_0)$. It follows that $K_x$ shares with $P$ the same basis of eigenvectors, so that $K_x\in \cK(P)$.
 Furthermore, we have that for any $l\in \lin \vert V\vert\rin$,
 \bq
 K_x[\varphi_l](x_0)&\df& 
 \sum _{z\in V}K_x(x_0,z)\varphi_l(z)\\&=& \sum _{z\in V}\sum_{k\in\lin \vert V\vert\rin} \varphi_k(x)\varphi_k(z)\varphi_l(z)\, \mu(z)\\
 &=&\varphi_l(x) \eq
 It implies that $K_x(x_0,\cdot)=\delta_x$.
 So for any $x\in V$, $\delta_x\in \cK(P,x_0)$. Taking into account that $\cK(P,x_0)$ is always a convex set, we get that $x_0\in\cH(P)$.\wwtbp
 \par\sm
 \begin{rem}
 (a) The uniplicity assumption cannot be removed for the reverse implication of Lemma~\ref{cHgroup}.
 Consider $P$ the transition kernel of the random walk on $V\df \ZZ/(n\ZZ)$, with $n\in\NN\setminus\{1,2\}$.
 At the end of Section 2.5 from \cite{MR2483738},
Bakry and Huet  show that $P$ does not satisfy the hypergroup property.
Nevertheless, consider for $v\in\ZZ/(n\ZZ)$, the translation by $v$ kernel $K$
defined by
\bq
\fo x,y\in\ZZ/(n\ZZ),\qquad K(x,y)&\df& \delta_{x+v}(y)\eq
Clearly $K\in\cK(P)$ and $K(0,\cdot)=\delta_v$, so that $\delta_v\in \cK(P,0)$ for all $v\in \ZZ/(n\ZZ)$.
It follows that $0\in\cH(P)$. More precisely, we have $\cH(P)=\ZZ/(n\ZZ)$.
\par\sm
(b) The example in (a) satisfies the \textbf{complex hypergroup property} with respect to any point $x_0\in \ZZ/(n\ZZ)$ (see Proposition 2.10 of Bakry and Huet \cite{MR2483738}), in the sense that
we can find an unitary basis $(\varphi_1,\ \varphi_2, ...,\ \varphi_{\lve V\rve})$ of $\LL^2(\mu,\CC)$
consisting of eigenvectors of $P$ 
  such that $\varphi_k(x_0)\not=0$ for all $k\in\lin\vert V\vert\rin$, and
\bqn{Bakry2}
\fo x,y,z\in V,\qquad \sum_{k\in\lin \vert V\vert\rin} \frac{\varphi_k(x)\varphi_k(y)\overline{\varphi_k(z)}}{\varphi_k(x_0)}&\geq &0\eqn
\par
So maybe the condition
\bqn{ghp}
\cH(P)&\not=&\emptyset\eqn is  related to the complex hypergroup property.
But here we will not investigate this question. We will mainly be interested in
\eqref{ghp}, seen as a generalization
 of the hypergroup property,
because it could   be considered for Markov kernels which are not reversible (or defined on abstract measurable spaces: \eqref{ghp} enables to avoid the technical difficulties
 related to the summations appearing in \eqref{Bakry} or \eqref{Bakry2} when the state space is not finite). 
 \end{rem}
 An irreducible birth and death kernel is necessarily uniplicit, so in the context of Theorem \ref{theo1},
 the hypergroup property for a Markov kernel $P$ is equivalent to \eqref{ghp}.
 We are thus lead to investigate the corresponding Markov commutator convex semi-group
 and will do it using general arguments.   The two properties 
 we will need are
 \begin{pro}
 \label{pro1}
 Assume that $P$ is an irreducible Markov kernel and let $\mu$ be its invariant probability.
 Then we have
 \bq
 \fo x\in \cH(P), \qquad \mu(x)&=&\min_V\mu\eq\
 \end{pro}\par
 For the second property, we need to introduce the \textbf{symmetry group} $\cS_P$ associated to $P$:
 it is the set of bijective mappings $g\st V\ri V$ such that
 \bqn{cSP}
 \fo x,y\in V,\qquad P(g(x),g(y))&=&P(x,y)\eqn
 For instance, one recovers the permutation group $\cS_{ V}$ of $V$  if $P$ is either the identity matrix $I$ (no move is permitted)
 or the matrix whose all off-diagonal entries are equal to $1/(\lve V\rve -1)$ (all ``true" moves are equally permitted). Indeed $\cS_P=\cS_{ V}$ if and only if
 $P$ is a convex combination of the two previous matrices, situations where all the elements of $V$ are indistinguishable with respect to the evolution dictated by $P$.\par
 \begin{pro}\label{pro2}
 Assume that $P$ is an uniplicit Markov kernel and let $x_0,x_1\in \cH(P)$.
 Then there exists $g\in\cS_P$ such that $g(x_1)=x_0$.
 Conversely, any $g\in \cS$ stabilizes $\cH(P)$, so that $\cH(P)$ is the orbit of any of its element under $\cS_P$.
 \end{pro}
 \par\sm
 Another natural question in the finite birth and death setting is the transposition of the Achour-Trimèche's theorem known  in the continuous framework.
 We did not succeed in getting a really satisfactory answer
 in this direction. The next result is obtained by adapting the arguments of Bakry and Huet \cite{MR2483738}.
 Let $\wi\cC$ be the subset of $U\in\cC$ such that $U(x+2)-U(x+1)\geq U(x+1)-U(x)+2\ln(2)$ for all $x\in\lin 0,N-2\rin$ (equivalently,
 $U$ is the restriction to $\lin 0,N\rin$ of a $\cC^2$ function on $[0,N]$ satisfying $U''\geq 2\ln(2)$).
 Let $\wi\cC_{\mathrm{m}}$  be the subset of $\wi \cC$ consisting of monotonous mappings such that 
 $\lve U(N)-U(N-1)\rve\wedge \lve U(1)-U(0)\rve \geq 2\ln(2)$. Consider also
 $\wi\cC_{\mathrm{s}}$ the subset of $\wi \cC$ consisting of
 mappings symmetric with respect to $N/2$.
\begin{pro}\label{pro3}
For any $U\in \wi\cC_{\mathrm{m}}\cup\wi\cC_{\mathrm{s}}$, the Metropolis kernel $M_U$ defined in \eqref{MU}
satisfies the hypergroup property. 
Thus the 
mapping $\wi\cC_{\mathrm{m}}\cup \wi\cC_{\mathrm{s}}\ni U\mapsto M_U$ is a birth and death 
Metropolis procedure
satisfying the hypergroup property. 
\end{pro}
In the one-dimensional  diffusive setting, the result corresponding to $\wi\cC_{\mathrm{m}}$ is due to Chebli \cite{MR0374545}.\par
Note that from Propositions \ref{pro1} and \ref{pro2}, we deduce that
 in the symmetric situation, $\cH(M_U)=\{0,N\}$, and that in the monotonous case with $U$ non-constant,
$\cH(M_U)$ is the singleton consisting of the boundary element with the smallest weight with respect to the reversible measure $\mu_U$.\par\sm
Remark \ref{finale} (d) 
gives another example of a generalized Metropolis procedure satisfying the hypergroup property for some convex potentials (more general
than those considered in Proposition~\ref{pro3}).
It would be very interesting to find other closed subsets $\cC'\subset \cC$ 
for which we can find a generalized Metropolis procedure $\cC'\ni U\mapsto P_U$ satisfying the hypergroup property (or to describe $\cC'\df\{U\in\cC\st \cH(M_U)\not=\emptyset\}$).
Especially to try to deduce the analogous results in the continuous framework, in order to recover Gasper's example \cite{MR0284628,MR0310536},
see also Bakry and Huet \cite{MR2483738} and Carlen, Geronimo and Loss \cite{MR2764893}.
\par\sm
From general considerations  related to the  Markov commutator convex semi-groups, we will also deduce the following criterion.
Let $\bar P$  be a Markov kernel on the finite set $\bar V$, 
 consider $\bar G$ a subgroup of $\cS_{\bar P}$ and denote by $\equiv$ the equivalence relation it induces on $\bar V$
via
\bq
\fo \bar x,\bar y\in\bar V,\qquad \bar x\equiv \bar y&\Leftrightarrow& \ex g\in\bar G\st g(\bar x)=\bar y\eq
Denote by $V$ the set of equivalence classes for $\equiv$ and by $\pi\st \bar V\ri V$ the associated projection mapping.
It is immediate to check that  a Markov kernel $P$ is well-defined on $V$ through the formula
\bq
\fo x,y\in V,\qquad P(x,y)&\df& \bar P(\bar x,\pi^{-1}(y))\eq
where $\bar x$ is any point of $\bar V$ such that $\pi(\bar x)=x$.
This construction corresponds to a reduction of the symmetries of $\bar P$. The next result shows that some properties of $\bar P$ are preserved under this operation.
It will be used to check the hypergroup property of $M_U$ for $U\in\wi\cC_{\mathrm{m}}$,
knowing it for $U\in\wi\cC_{\mathrm{s}}$.
\begin{pro}\label{pro4}
Assume that $\bar P$  is uniplicit  and  satisfies Condition \eqref{ghp}.
Then the same remains true for $P$.
\end{pro}
\par
If the uniplicity of $\bar P$ could be removed from this statement and be replaced by the uniplicity of $P$  (this is a weaker condition, since it will be seen in the proof of Corollary \ref{cor1}
 that the uniplicity of $\bar P$ implies that of $P$ under the assumptions of Proposition \ref{pro4}),
this result would provide an abstract rewriting in the finite context
 of the Carlen, Geronimo and Loss method  \cite{MR2764893}.
 This conjectured extension seems quite challenging, 
 some assumptions could be required on the subgroup $G$.
Maybe they do not appear here, because when $\bar P$ is uniplicit, $\cS_{\bar P}$ is commutative, see Remark \ref{commut} (a) below.
 \par\me
 In the next section we will study   the 
 Markov commutator convex semi-group in the general finite framework, obtaining in particular  Propositions \ref{pro1},  \ref{pro2} and \ref{pro4}.
Advantage will be taken of the relations between the 
 Markov commutator convex semi-group and the theory of Markov intertwining as it was
 developed by Diaconis and Fill \cite{MR1071805}.
 In the last section we consider more specifically the birth and death case and prove Theorem \ref{theo1} and Proposition \ref{pro3}.

\section{General properties}

This is the beginning of a systematic investigation of the Markov commutator convex semigroup $\cK(P)$ associated to a finite Markov kernel $P$.\par\me
We start by recalling some elements of the theory of Markov intertwining due to Diaconis and Fill \cite{MR1071805}.
Let $X\df (X_n)_{n\in\ZZ_+}$ and $\bar X\df (\bar X_n)_{n\in\ZZ_+}$ be two Markov chains, respectively on the finite state spaces
$V$ and $\bar V$. The respective transition kernels are denoted $P$ and $\bar P$, and the initial distributions $m_0$
and $\bar m_0$.
We say that $X$ is \textbf{intertwined} with $\bar X$ \textbf{through the Markov link} $\Lambda$ (which is a Markov kernel from $\bar V$ to $V$,
seen as a $\bar V\times V$ matrix), if there is a coupling $(X,\bar X)$ such that the two following conditions are met:
\bqn{d1}
\fo n\in\ZZ_+,\qquad \cL(\bar X_{\lin 0,n\rin}\vert X)&=& \cL(\bar X_{\lin 0,n\rin}\vert X_{\lin 0,n\rin})
\eqn
where as usual this identity of conditional laws has to be understood a.s.\ with respect to the probability measure underlying the coupling. The trajectorial notation 
$X_{\lin 0,n\rin}\df (X_p)_{p\in\lin 0,n\rin}$ was used.
\bqn{d2}
\fo n\in\ZZ_+,\qquad \cL( X_{n}\vert \bar X_{\lin 0,n\rin})&=& \Lambda(\bar X_n,\cdot)\eqn
When these assumptions are satisfied, we write $X\prec_\Lambda \bar X$ and $\bar X$ is also said to be a \textbf{dual chain} of $X$ through $\Lambda$.
The notation $X\prec \bar X$ will notify there exists $\Lambda$ such that $X\prec_\Lambda \bar X$.\par
We say that $(m_0,P)$ is \textbf{intertwined} with  $(\bar m_0,\bar P)$ \textbf{through the Markov link} $\Lambda$ 
if 
\bqn{equi}
m_0\ =\ \bar m_0\Lambda&\hbox{and}& \bar P\Lambda\ =\ \Lambda P\eqn
We denote this relation by $( m_0,P)\prec_\Lambda (\bar m_0,\bar P)$ and as above, $( m_0,P)\prec (\bar m_0,\bar P)$ means there exists a kernel $\Lambda$
such that \eqref{equi} is satisfied. 
\par
Diaconis and Fill \cite{MR1071805} have shown that these notions of  intertwining coincide, at least if $\bar X$ visits the whole state space $\bar V$ (in particular if $\bar P$ is irreducible):
\begin{pro}\label{DF}
With the above notations, we have
\bq
( m_0,P)\prec_\Lambda (\bar m_0,\bar P)&\Rightarrow& X\prec_\Lambda \bar X
\eq
Furthermore if for any $\bar x\in\bar V$, there exists $n\in\ZZ_+$ such that $\PP[\bar X_n= \bar x]>0$, then
\bq
X\prec_\Lambda \bar X&\Rightarrow& ( m_0,P)\prec_\Lambda (\bar m_0,\bar P)
\eq
\end{pro}
\par
\prooff{Proof:}
More specifically, the construction of the coupling of $X$ and $\bar X$ satisfying the conditions \eqref{d1} and \eqref{d2} under the assumption $( m_0,P)\prec_\Lambda (\bar m_0,\bar P)$
is described in Theorem 2.17 of Diaconis and Fill \cite{MR1071805}. The other implication can also be deduced from their considerations. For the sake of completeness, here are
some arguments, directly based on  the hypotheses \eqref{d1} and \eqref{d2}.
\par
From \eqref{d2}, we deduce that for all $n\in\ZZ_+$, $\cL( X_{n}\vert \bar X_n)= \Lambda(\bar X_n,\cdot)$ so that by integration with respect to $\bar X_n$,
we get 
$\cL(X_n)=\cL(\bar X_n)\Lambda$. In particular for $n=0$, we obtain $m_0 = \bar m_0\Lambda$.
\\
Let $f$ and $\bar f$ two test functions defined respectively on $V$ and $\bar V$.
For fixed $n\in\ZZ_+$, we compute $\EE[\bar f(\bar X_{n}) f(X_{n+1})]$ in two ways.
First, using \eqref{d2} and the Markov property of $\bar X$,
\bq
\EE[\bar f(\bar X_{n}) f(X_{n+1})]&=&\EE[\bar f(\bar X_{n}) \EE[f(X_{n+1})\vert \bar X_{\lin 0, n+1\rin}]]\\
&=&\EE[\bar f(\bar X_{n}) \Lambda[f](\bar X_{n+1})]\\
&=&\EE[\bar f(\bar X_{n})(\bar P\Lambda)[f](\bar X_{n})]\eq
Second, 
using \eqref{d1} and the Markov property of $X$,
\bq
\EE[\bar f(\bar X_{n}) f(X_{n+1})]&=&\EE[\EE[\bar f(\bar X_{n}) \vert X]f(X_{n+1})]\\
&=&\EE[\EE[\bar f(\bar X_{n}) \vert X_{\lin 0,n\rin}]f(X_{n+1})]\\
&=&\EE[\EE[\bar f(\bar X_{n}) \vert X_{\lin 0,n\rin}]P[f](X_{n})]\\
&=&\EE[\bar f(\bar X_{n})P[f](X_{n})]\\
&=&\EE[\bar f(\bar X_{n})\EE[P[f](X_{n})\vert \bar X_{\lin 0,n\rin}]]\\
&=& \EE[\bar f(\bar X_{n})(\Lambda P)[f](\bar X_{n})]\eq
Since this is true for any $\bar f$, we deduce that a.s., 
\bq
(\Lambda P)[f](\bar X_{n})&=&(\bar P\Lambda)[f](\bar X_{n})\eq
and due to the assumption on $\bar X$,
\bq
\fo \bar x\in \bar V,\qquad 
(\Lambda P)[f](\bar x)&=&(\bar P\Lambda)[f](\bar x)\eq
Since it is true for all $f$, it follows that $\Lambda P=\bar P\Lambda$.\wwtbp
\begin{rem}\label{relation}
(a) The relation $\prec$ is clearly reflexive (through the identity link) and it can be easily checked to be transitive (for instance at the level of the Markov chains,
if $X\prec_\Lambda X'$ and $X'\prec_{\Lambda'} X''$ then  $X\prec_{\Lambda\Lambda'} X''$).
Thus $\prec$ is a pre-order,
e.g.\ on the trajectorial laws of finite Markov chains (whose state space is a subset of $\NN$, to work on a defined set).
It is then tempting to verify if it would not be  an equivalence  or an order relation.
To see that $\prec$ is none, consider $Y$ the trivial Markov chain on a singleton.
For any finite Markov chain $X$, we have $Y\prec X$, but $X\prec Y$ is equivalent to the stationarity of $X$ (namely the initial
distribution of $X$ is invariant for its transition kernel).
It follows that $\prec$ is neither symmetrical nor anti-symmetrical.
Next, one can define an equivalence relation $X\sim X'$ via $X\prec X'$ and $X'\prec X$.
On the corresponding equivalence classes, $\prec$ defines a partial order relation,
in some sense it should compare the difficulty of reaching an equilibrium (see also Remark \ref{sepdis} below). The ``stationarity" class of the trivial chain $Y$ is minimal for this order.
\par\sm
(b) Similar conditions are valid for the algebraic intertwining between couples consisting of a probability measure and a Markov kernel.
If the finite state set $V$ and the Markov kernel $P$ are fixed, we induce a relation on $\cP(V)$ via $m_0\prec \bar m_0$ if and only if $(m_0, P)\prec (\bar m_0,P)$.
It can be transformed into 
 an order relation on $\cP(V)/\sim$ by introducing an equivalence relation $\sim$ as above. It  heuristically corresponds to the proximity to the set of invariant measures for $P$,
 which are the minimal elements. Note that the semigroup $(P^n)_{n\in\ZZ_+}$ is non-increasing with respect to $\prec$, since we have
 $(m_0P,P)\prec_P (m_0,P)$.
\end{rem}
The main interest of associating a dual chain $\bar X$ to a given Markov chain $X$ is
that it enables to construct strong times (see for instance Diaconis and Fill \cite{MR1071805}, Fill \cite{MR2530102},
Diaconis and Miclo \cite{MR2530103} and 
\cite{MR2654550}).
A stopping time $\tau$ for $X$ (with respect to a filtration containing the filtration generated by $X$)
is a \textbf{strong time} if it is a.s.\ finite and if $\tau$ and $X_\tau$ are independent. The basic principle of the construction is
the following well-known result, whose proof is  given for the sake of completeness.
\begin{lem}\label{st}
Let $(X,\bar X)$ be a coupling satisfying  \eqref{d2}, then this equality can be extended to any a.s.\ finite stopping time $\tau$ for $\bar X$,
namely
\bq
\cL( X_{\tau}\vert \bar X_{\lin 0,\tau\rin})&=& \Lambda(\bar X_\tau,\cdot)\eq
If in addition $(X,\bar X)$  satisfies  \eqref{d1}, then $\tau$ is a strong time if  $\Lambda(\bar X_\tau,\cdot)$ is independent from $\tau$
(for instance if $\Lambda(\bar X_\tau,\cdot)$ ``is not really depending on" $\bar X_\tau$, e.g.\ if $\bar X_\tau$ is a.s.\ equal to a fixed point).
\end{lem}
\proof
The first assertion is an outcome of the notion of a stopping time:
Let $f$ be a function defined on $V$ and $\bar F$ a bounded functional measurable with respect to the stopped trajectory
$ \bar X_{\lin 0,\tau\rin}$. We compute that
\bq
\EE[f(X_\tau)\bar F]&=&\sum_{n\in\ZZ_+}\EE[f(X_n)\bar F\un_{\tau=n}]\\
&=&\sum_{n\in\ZZ_+}\EE[\EE[f(X_n)\vert \bar X_{\lin 0,n\rin}]\bar F\un_{\tau=n}]\\
&=&\sum_{n\in\ZZ_+}\EE[\Lambda[f](\bar X_n)\bar F\un_{\tau=n}]\\
&=&\EE[\Lambda [f](\bar X_\tau)\bar F]\eq
where the second equality comes from the fact that $\bar F\un_{\tau=n}$ is measurable with respect to $\bar X_{\lin 0,n\rin}$.
The first wanted result follows, since this is true for all $f$ and $\bar F$ as above.\par
For the second assertion, note that \eqref{d1} implies that a stopping time for $\bar X$ is also a stopping time for $X$.
Let  $f$ be a function defined on $V$ and let $g$ be a bounded measurable mapping on $\RR_+$.
Since $\tau$ is measurable with respect to $ \bar X_{\lin 0,\tau\rin}$, we have
\bq
\EE[f(X_\tau)g(\tau)]&=&\EE[\EE[f(X_\tau)\vert \bar X_{\lin 0,\tau\rin}] g(\tau)]\\
&=& \EE[\Lambda[f](\bar X_\tau)g(\tau)]\\&=&\EE[\Lambda[f](\bar X_\tau)]\EE[g(\tau)]\\
&=&\EE[f(X_\tau)]\EE[g(\tau)]\eq
where the third equality comes from the assumption made on $\Lambda(\bar X_\tau,\cdot)$.
The independence of $\tau$ and $X_\tau$ follows, since $f$ and $g$ were arbitrary.\wwtbp\par
For the purpose of proving Proposition \ref{pro1}, we will only use the first part of the above lemma, 
even if the stopping times we will consider are indeed strong times.\par
\sm
Indeed, it is time to come back to the
Markov commutator convex semigroup $\cK(P)$ associated to an irreducible finite Markov kernel $P$.
Denote $X^{m_0}\df (X_t^{m_0})_{t\geq 0}$ a Markov chain with $P$ as transition kernel and $m_0\in \cP(V)$ for initial law.
From the definitions and Proposition \ref{DF}, we have for any $K\in \cK(P)$ and any initial distribution $m_0$, $X^{m_0K}\prec_K X^{m_0}$.
Thus it appears that $x_0\in\cH(P)$ if and only if for any $m_0\in\cP(V)$, there exists a Markov kernel $K$ on $V$ such that
$X^{m_0}\prec_K X^{x_0}$ (as customary, $X^{x_0}$ is a shorthand for $X^{\delta_{x_0}}$). 
In particular, if $P$ is uniplicit, then $P$ satisfies the hypergroup property if and only there exists $x_0\in V$ such that
 for any $m_0\in\cP(V)$, $X^{m_0}\prec X^{x_0}$.
 More generally, we get the following interpretation:
 \bq\fo x\in V,\qquad \cK(P,x)&=&\{m\in \cP(V)\st X^m\prec X^x\}\eq
  \par
All preliminaries are now in place for the
\prooff{Proof of Proposition \ref{pro1}}
Consider $x_0\in\cH(P)$ and let $x_1$ be any point of $V$. We want to show that 
$\mu(x_0)\leq \mu(x_1)$. \\
By definition of $\cH(P)$, there exists $K\in\cK(P)$ such that $K(x_0,\cdot)=\delta_{x_1}$, so that from Proposition \ref{DF},
$X^{x_1}\prec_KX^{x_0}$, i.e.\ 
we can construct a coupling of $X^{x_0}$ and $X^{x_1}$
satisfying \eqref{d1} and \eqref{d2} with $\Lambda\df K$.
\\
Let $(\tau_n)_{n\in\ZZ_+}$ be the sequence of stopping times for $X^{x_0}$ defined by iteration through
$\tau_0=0$ and
\bq
\fo n\in\ZZ_+,\qquad \tau_{n+1}&\df& \inf\{p>\tau_n\st X_p=x_0\}\eq
According to Lemma \ref{st}, for any $n\in\ZZ_+$,
\bq
\cL(X^{x_1}_{\tau_n}\vert X^{x_0}_{\lin 0,\tau_n\rin})&=& \delta_{x_1}\eq
It means that each time $X^{x_0}$ is in $x_0$, then $X^{x_1}$ is in $x_1$.
It remains to apply the ergodic theorem to get
\bq
\mu(x_0)&=&\lim_{n\ri\iy}\frac1{n+1}\sum_{p\in\lin 0, n\rin}\un_{x_0}(X^{x_0})\\
&\leq & \lim_{n\ri\iy}\frac1{n+1}\sum_{p\in\lin 0, n\rin}\un_{x_1}(X^{x_1})\\
&=&\mu(x_1)\eq
where the (in)equalities are valid a.s.\wwtbp
The elements of $\cH(P)$ satisfies other optimization properties, they are for instance points
from which it is the most difficult to reach equilibrium in the separation discrepancy sense:
\begin{rem}\label{sepdis}
Recall that the \textbf{separation discrepancy} $s(m,\mu)$ between two probability measures on $V$ is defined by
\bq
s(m,\mu)&\df& \sup_{x\in V}1-\frac{m(x)}{\mu(x)}\eq
(with the usual convention: $r/0=+\iy$ for any $r>0$, but $0/0=0$).
\\
A \textbf{stationary time} $\tau$ for an irreducible Markov chain $X^{m_0}\df(X^{m_0}_n)_{n\in\ZZ_+}$ ($m_0$ still stands  for the initial distribution)  is a strong time such that $X^{m_0}_\tau$ is distributed according to the associated invariant measure $\mu$.
Aldous and Diaconis  \cite{MR876954} have shown that if the transition kernel is aperiodic and irreducible, then for any initial distribution $m_0$,
there exists a stationary time $\tau^{m_0}$ associated to $X^{m_0}$ satisfying
\bq
\fo n\in\ZZ_+,\qquad \PP[\tau^{m_0}> n]&=&s(m_0P^n,\mu)\eq
Furthermore $\tau^{m_0}$ is stochastically smaller than any stationary time associated to $X^{m_0}$.
\par
The proof of Proposition \ref{pro1} can be slightly modified to show that if $x_0\in\cH(P)$, then
$\tau^{x_0}$ is stochastically larger than $\tau^{m_0}$ for any initial distribution $m_0$.
Indeed, if $K\in\cK(P)$ is such that $K(x_0,\cdot)=m_0$, then considering a coupling of $X^{x_0}$ and $X^{m_0}$
realizing the relation $X^{m_0}\prec_K X^{x_0}$, it appears that $\tau^{x^0}$ is a stationary time for $X^{m_0}$.
It is a consequence of the fact that all the elements of $\cK(P)$  admit $\mu$ for invariant measure, as it was seen in \cite{miclo:hal-01117051}
(only the irreducibility of $P$ is needed for this property).
The stochastic domination of $\tau^{m_0}$  by $\tau^{x_0}$ ensures that for any initial distribution $m_0$
(or equivalently for any Dirac mass $m_0=\delta_{x_1}$, with $x_1$ in the state space $V$),
\bq
\fo n\in\ZZ_+,\qquad s(m_0P^n,\mu)&\leq & s(P^n(x_0,\cdot),\mu)\eq
\end{rem}
\par
\sm
To go in the direction of Proposition \ref{pro2}, we begin by a simple technical result:
\begin{lem}\label{inversion}
Let $K$ and $K'$ be two Markov kernels on $V$ such that $K'K=I$, the identity kernel. Then there exist $g\in \cS_V$
such that
\bq
\fo x,y\in V,\qquad \lt\{\begin{array}{rcl}
K(x,y)&=&\delta_{g(x)}(y)\\
K'(x,y)&=&\delta_{g^{-1}(x)}(y)\end{array}\rt.\eq
\end{lem}
\proof
By contradiction, assume  there exists $x\in V$ such that $K(x,\cdot)$ is not a Dirac mass.
Then for any $y\in V$, if $K'(y,x)>0$ then $K'K(y,\cdot)$ cannot be a Dirac mass.
This is not compatible with $K'K=I$, so we must have $K'(y,x)=0$ for all $y\in V$.
It implies that $K'$ is not invertible, in contradiction again with our assumption.
So for any $x\in V$, $K(x,\cdot)$ is a Dirac mass $\delta_{g(x)}$ for some $g(x)\in V$.
Since $K$ is invertible, necessarily the mapping $g$ is also invertible.
The announced result follows at once.\wwtbp
\par
In addition, we will need the following consequence of the uniplicit assumption.
\begin{lem}\label{unic}
Assume that $P$ is uniplicit, then for any fixed $x_0\in \cH(P)$,  the affine mapping
\bq
\cK(P)\ni K&\mapsto & K(x_0,\cdot)\in \cK(P,x_0)\eq
is one-to-one. 
\end{lem}
\proof
Fix $x_0\in \cH(P)$ and $m_0\in\cP(V)$, it is sufficient to see there is exactly one matrix $K$ solution
to the equations
\bq
K(x_0,\cdot)&=&m_0\\
KP&=&PK\eq
Indeed, consider $\mu$ the reversible probability for $P$ and let $\varphi_1,\ \varphi_2, ...,\ \varphi_{\lve V\rve}$ be an orthonormal (in $\LL2(\mu)$) basis of  eigenvectors 
associated to $P$  as in the introduction.
By the commutation of $K$ with $P$, this is also a basis of  eigenvectors for $K$. Thus we can find numbers $a_1, a_2, ..., a_{\lve V\rve}$ such that
\bq
\fo x,y\in V,\qquad K(x,y)&=&\sum_{l\in\lin 1,\lve V\rve\rin} a_l\varphi_l(x)\varphi_l(y)\mu(y)\eq
The first condition then reads
\bq
\fo y\in V,\qquad \frac{m_0}{\mu}(y)&=&\sum_{l\in\lin  1,\lve V\rve\rin} a_l\varphi_l(x_0)\varphi_l(y)\eq
namely $( a_l\varphi_l(x_0))_{l\in \lin  1,\lve V\rve\rin}$ are the coefficients of $m_0/\mu$ in the basis
$(\varphi_1,\ \varphi_2, ...,\ \varphi_{\lve V\rve})$.
Since $\varphi_l(x_0)\not=0$ for all $l\in \lin  1,\lve V\rve\rin$, according to Lemma \ref{cHgroup},
we get that the $a_1, a_2, ..., a_{\lve V\rve}$ are uniquely determined.\wwtbp
In particular if $P$ is an uniplicit kernel satisfying the hypergroup property, then $\cK(P)$ is a simplex.
It is sometimes possible to go further:
\begin{rem} In fact 
the above proof shows that if $x_0\in V$ is any point such that
$\varphi_l(x_0)\not=0$ for all $l\in \lin  1,\lve V\rve\rin$, then
the conclusion of Lemma \ref{unic} still holds if $P$ is uniplicit.
If furthermore \eqref{ghp} holds, then $\cK(P)$  is a simplex
as well as each of the $\cK(P,x)$, for $x\in V$.\par
 Let $\cR$ be the set of Markov kernels which are irreducible and reversible.
It can be easily seen that the subset of elements of $\cR$ which are uniplicit and whose eigenvectors never vanish
 is a dense open subset of $\cR$. But since $\cH$, the subset of $\cR$ consisting of kernels satisfying the hypergroup
 property, is very slim in $\cR$, it is no longer clear whether or not the subset of elements of $\cH$ which are uniplicit and whose eigenvectors never vanish
 is  a dense open subset of $\cH$. If it was true, it could be concluded that ``generically", $\cK(P,x)$ is a simplex  for  $P\in\cH$ and  $x\in V$ .
\end{rem}\par
We  have all the ingredients for the
\prooff{Proof of Proposition \ref{pro2}}
Let be given $x_0,x_1\in \cH(P)$. Then 
there exist $K',K\in\cK(P)$ such that
\bqn{xx}
\nonumber K'(x_0,\cdot)&=&\delta_{x_1}\\
K(x_1,\cdot)&=&\delta_{x_0}\eqn
Thus we get that $K'K(x_0,\cdot)=\delta_{x_0}$.
Since $P$ is assumed to be uniplicit, we get from Lemma \ref{unic} that $K'K\in\cK(P)$ is uniquely determined by this relation.
It appears there is no alternative: $K'K=I$.
Lemma \ref{inversion} enables to find a permutation $g\in\cS_V$ such $K$ is the Markov kernel induced by $g$. 
Note that \eqref{xx} translates into $g(x_1)=x_0$.
The commutation of $K$ and $P$ then implies that
\bq
\fo x,y\in V,\qquad P(g(x),y)&=&P(x,g^{-1}(y))\eq
which can be rewritten under the form \eqref{cSP}
namely $g\in\cS_P$. The remaining assertions of Proposition~\ref{pro2} are straightforward.\wwtbp
\par\sm
We are  now 
going in the direction of Proposition \ref{pro4} through a sequence of general arguments, in the hope 
they present in a clear way the problems one will encounter in trying to generalize it.
We start by recalling some  considerations from  \cite{miclo:hal-01117051}.
A Markov kernel $\Lambda$ from $\bar V$ to $V$
can be interpreted as an operator sending any function $f$ defined on $V$ to the mapping $\Lambda[f]$ defined on $\bar V$ by
\bq
\fo \bar x\in\bar V,\qquad \Lambda[f](\bar x)&\df& \sum_{x\in V}\Lambda(\bar x,x)f(x)\eq
Let $\bar\mu$ be a probability measure given on $\bar V$ and consider $\mu\df \bar\mu\Lambda$ its image by $\Lambda$.
Then $\Lambda$ can be seen as an operator from $\LL^2(\mu)$ to $\LL^2(\bar\mu)$ (because $\Lambda[f]$ is $\bar\mu$-negligible
if $f$ is $\mu$-negligible). It enables to define $\Lambda^*$ its dual operator  from $\LL^2(\bar\mu)$ to $\LL^2(\mu)$,
which is Markovian in the sense that
\bq
\Lambda^*[\un_{\bar V}]&=&\un_{V}\\
\fo f\in\LL^2(\bar \mu), \quad f\geq 0&\Rightarrow& \Lambda^*[f]\geq 0\eq
where the relations have to be understood $\bar\mu$- or $\mu$-a.s.\par
If $\bar\mu$ and $\mu$ give positive weights to all points of $\bar V$ and $V$ respectively,
then $\Lambda^*$ can be seen as a Markov kernel from $V$ to $\bar V$.
\begin{rem}\label{dual}
In the intertwining framework, similar considerations are valid for $\bar P$ and $P$, in order to define $\bar P^*$ and $P^*$,
seen as Markov operators on  $\LL^2(\bar\mu)$ and $\LL^2(\mu)$,
when $\bar\mu$ and $\mu$ are invariant probability measures, respectively for $\bar P$  and $P$, i.e.\ $\bar\mu\bar P=\bar\mu$ and $\mu P= \mu$.
Thus to be able to consider $\bar P^*$ and $P^*$ as Markov matrices, it is convenient to make the following assumption: we say that the couple
$(\bar P,\Lambda)$ is \textbf{positive}, if $\bar P$ admits a positive invariant measure $\bar \mu$
and if $\mu\df \bar\mu\Lambda$ is also positive.
Up to reducing $\bar V$ and $V$ respectively to the support of $\bar\mu$ and $\mu$,
it is always possible to come back to this case.
 Note that  the commutation relation \bqn{intertw}
 \bar P\Lambda&=&\Lambda P\eqn
implies that $\mu$ is an invariant probability for $P$. \par
Under the hypotheses that $(\bar P,\Lambda)$ is positive and that \eqref{intertw}
is satisfied, we get a dual commutation relation:
\bq
P^*\Lambda^*&=&\Lambda^*\bar P^*\eq
If furthermore we assume that $(m_0,P)\prec_\Lambda (\bar m_0,\bar P)$ and that
\bqn{m0hyp}
\bar m_0\Lambda\Lambda^*&=&\bar m_0\eqn
then we get the intertwining relation
\bq
(m_0,\bar P^*)&\prec_{\Lambda^*} &(\bar m_0, P^*)\eq
The reversibility assumption for $\bar P$ with respect to $\bar\mu$ amounts to 
$\bar P^*=\bar P$ and similarly for $P$.
These considerations lead to a restricted symmetry property for the relation $\prec$:
 $(m_0,P)\prec_\Lambda (\bar m_0,\bar P)$ implies $(\bar m_0,\bar P)\prec_{\Lambda^*} (m_0, P)$
 under the assumptions that $(\bar P,\Lambda)$ is positive, that $\bar P$ and $P$ are reversible
 and that \eqref{m0hyp} is satisfied. This is an instance of the equivalence relation $\sim$ introduced in Remark \ref{relation}.\par
 We give below in Remark \ref{nondet} (b) a natural condition under which  \eqref{m0hyp} is true.
\end{rem}
Beyond reversibility or uniplicity, an important assumption will be 
\bqn{hyp}
\Lambda\Lambda^*\bar P\Lambda&=&\bar P\Lambda\eqn
(this condition for the Markov kernel $\bar P$ is an analogue of \eqref{m0hyp} for the probability measure $\bar m_0$).
Define \bqn{hyp2}
P&\df& \Lambda^* \bar P\Lambda\eqn
From \eqref{hyp},
it appears that $\bar P$ and $P$ are intertwined through $\Lambda$,
namely \eqref{intertw} is satisfied.
We can go further in the exploration of $\cK(P)$ with the help of $\cK(\bar P)$: the next result is a slight modification of Proposition~3 of \cite{miclo:hal-01117051}, where $\cK(\bar P)$
was replaced by the 
 smaller set
 \bq
 \cK(\bar P,\Lambda)&\df& \{K\in \cK(\bar P)\st \Lambda\Lambda^*\bar K\Lambda=\bar K\Lambda\}\eq
namely the set of elements from $\cK(\bar P)$  satisfying the condition \eqref{hyp}.
It is also a convex semigroup and in Lemma \ref{cKPL} some conditions will be given so that $ \cK(\bar P,\Lambda)=\cK(\bar P)$.
\begin{lem}\label{CGL}
Assume that $\bar P$ is reversible with respect to $\bar\mu$ and that \eqref{hyp} holds, then we have 
\bq
 \Lambda^*\cK(\bar P)\Lambda &\subset & \cK(P)\eq
 \end{lem}
 \proof
 For any $\bar K\in\cK(\bar P)$, we compute that
 \bqn{KK}
\nonumber \Lambda^* \bar K \Lambda P&=& \Lambda^* \bar K \bar P \Lambda\\
 \nonumber &=& \Lambda^*  \bar P \bar K\Lambda\\
\nonumber  &=& \Lambda^*  \bar P\Lambda \Lambda^* \bar K\Lambda\\
&=&P  \Lambda^* \bar K\Lambda\eqn
where for the third equality, we have used the dual relation of \eqref{hyp}
asserting that 
$\Lambda^*\bar P^*\Lambda\Lambda^*=\Lambda^*\bar P^*$,
namely
$\Lambda^*\bar P\Lambda\Lambda^*=\Lambda^*\bar P$,
since $\bar P=\bar P^*$.
Relation \eqref{KK} shows that $\Lambda^* \bar K\Lambda$ belongs to $\cK(P)$.\wwtbp
\par
Condition \eqref{hyp} seems quite strange at first view and we would have liked to  only work with \eqref{intertw}.
Lemma~\ref{hyphyp} below will show this is possible when $\Lambda$ is deterministic.
\par
\sm
The motivation for Proposition 3 of \cite{miclo:hal-01117051}
was to give an abstract version  in the finite context of a method of
Carlen, Geronimo and Loss \cite{MR2764893} to recover the hypergroup property in the context of Jacobi polynomials,
result initially due to Gasper \cite{MR0284628,MR0310536}. The underlying idea is equally conveyed by Lemma \ref{CGL}:
to prove \eqref{ghp}, one tries to find a Markov model (or several ones) $\bar P$, above $P$ in the sense of intertwining
(namely according to the order relation induced by $\prec$ as in Remark \ref{relation}), such that
 $\cK(\bar P)$ 
is relatively easy to apprehend. If it appears  that $\cK(\bar P)$
is
 quite big, then the inclusion of Lemma  \ref{CGL} gives an opportunity to show that $\cK(P)$ is also big, leading us toward   \eqref{ghp}.
 But to guess such a nice Markov kernel $\bar P$ from $P$ may not be  an easy task!
That is why we now go in the reverse direction, starting with $\bar P$.
In particular it is natural to wonder when does 
 \bqn{bghp}
\cH(\bar P)&\not=& \emptyset\eqn
imply  \eqref{ghp}. 
Before partially answering this question, let us mention a construction of Markov kernels satisfying \eqref{bghp}.
\begin{rem}\label{tensorization}
(a) Any irreducible Markov kernel $P$ on $\{0,1\}$ satisfies \eqref{ghp}. Indeed, let $\mu$ be the associated invariant measure
and by symmetry, assume that $\mu(0)\leq \mu(1)$.
Then there exists $a\in[-\mu(0)/\mu(1),1]$ such that $P=aI+(1-a)\mu$, where $\mu$ is seen as the Markov kernel 
whose two rows are equal to $\mu$.
Any Markov kernel $K\df bI+(1-b)\mu$, with $b\in[-\mu(0)/\mu(1),1]$, belongs to $\cK(P)$.
Taking $b=-\mu(0)/\mu(1)$ (respectively $b=1$), the first row of $K$ is $(0,1)$ (resp.\ $(1,0)$).
This shows that $0\in \cH(P)$.\par\sm
(b) If $P_1$ and $P_2$ are two Markov kernels on $V_1$ and $V_2$, then $P_1\otimes P_2$ is a Markov kernel
on $V_1\times V_2$. It appears that $\cK(P_1)\otimes\cK(P_2)\subset \cK(P_1\otimes P_2)$
and in particular $\cH(P_1)\times \cH(P_2)\subset \cH(P_1\otimes P_2)$.
\par\sm
(c) From the two points above, it follows that if $P$ is an irreducible Markov kernel on $\{0,1\}$, then for any $N\in\NN$,
$\bar P\df P^{\otimes N}$ satisfies \eqref{bghp}. Such Markov kernels were used in \cite{miclo:hal-01117051}
to recover the hypergroup property of the biased Ehrenfest model (initially due to Eagleson \cite{MR0328162}).
\end{rem}
We introduce now three assumptions which are helpful in the direction of deducing \eqref{ghp} from \eqref{bghp}.
 \par
 First, the surjectivity of $\Lambda$ as an operator on $\cP(\bar V)$:
 \bqn{H1}
\cP(\bar V)\Lambda&=&\cP(V)
\eqn
Second, the determinism of $\Lambda$ on $\cH(\bar P)$:
 \bqn{H2}
\fo \bar x_0\in\cH(\bar P),\qquad 
\Lambda(\bar x_0,\cdot)&=&\delta_{\pi(\bar x_0)}
\eqn
where $\pi(\bar x_0)$ is an element of $V$. Denote $\pi(\cH(\bar P))$ the image by $\pi$ of $ \cH(\bar P)$.
The last hypothesis is an extension of \eqref{hyp} to the identity kernel:
\bqn{H3}
\Lambda\Lambda^*\Lambda&=&\Lambda\eqn
Note that by multiplication on the left or on the right by $\Lambda^*$, this implies that $\Lambda\Lambda^*$ and $\Lambda^*\Lambda$ are projection operators in their respective spaces  $\LL^2(\bar\mu)$ and $\LL^2(\mu)$.
\begin{pro}\label{trop}
Assume $\bar P$ is uniplicit  and 
\eqref{hyp},  \eqref{bghp}, \eqref{H1}, \eqref{H2} and \eqref{H3} hold. Then
 \eqref{ghp} is satisfied with $P$ given by \eqref{hyp2}
and more precisely $\pi(\cH(\bar P))\subset \cH(P)$.
\end{pro}
\par
Before proving this statement, let us give another important consequence of uniplicity.
If $P$ is a Markov kernel on $V$, let $\cA(P)$ be the algebra generated by $P$,
namely the set of finite combinations of the
 form
 $a_0I+a_1P+a_2P^2+\cdots+ a_nP^n$, where $n\in\ZZ_+$ and $a_0, a_1, a_2, ..., a_n\in \RR$.
 Denote also by $\cK(V)$ the convex set of Markov kernels on $V$.
 \begin{lem}\label{cKPL}
 Assume that $P$ is uniplicit. Then we have
 \bq
 \cK(P)&=&\cA(P)\cap\cK(V)\eq
In particular if \eqref{H3} holds and $\bar P$ is uniplicit and satisfies \eqref{hyp},
then the latter property can be extended to $\cK(\bar P)$:
\bqn{hypK}
\fo \bar K\in\cK(\bar P),\qquad \Lambda\Lambda^* \bar K\Lambda&=&\bar K\Lambda\eqn
\end{lem}
\proof
Let  $(\varphi_1,\ \varphi_2, ...,\ \varphi_{\lve V\rve})$ be an
 orthonormal basis of  eigenvectors of $P$ and let $\lambda_1,\ \lambda_2, ...,\ \lambda_{\lve V\rve}$
 be the corresponding eigenvalues.
 Consider $K\in\cK(P)$, by commutativity, $(\varphi_1,\ \varphi_2, ...,\ \varphi_{\lve V\rve})$ is also
 a basis of eigenvectors of $K$, denote by $\theta_1,\ \theta_2, ...,\ \theta_{\lve V\rve}$
 the associated eigenvalues. Since the $\lambda_1,\ \lambda_2, ...,\ \lambda_{\lve V\rve}$
 are all distinct, we can find a polynomial $R$ of degree at most $\lve V\rve$ such that
 \bq
 \fo l\in\lin \lve V\rve\rin,\qquad R(\lambda_l)&=&\theta_l\eq
 It follows that $K=R(P)$, showing that $\cK(P)\subset \cA(P)\cap \cK(V)$.
 The reverse inclusion is obviously always true.\par
 The second assertion of the lemma comes from the fact that \eqref{hyp} implies that 
 \bq
 \fo n\in\NN,\qquad 
   \Lambda\Lambda^* \bar P^n\Lambda&=& \bar P^n\Lambda\eq
Indeed, this is shown by iteration on $n\in\NN$:
\bq
 \Lambda\Lambda^* \bar P^{n+1}\Lambda&=& \Lambda\Lambda^* \bar P^{n}(\bar P\Lambda)\\
 &=&
 \Lambda\Lambda^*   \bar P^n (\Lambda\Lambda^*\bar P\Lambda)\\
 &=&
( \Lambda\Lambda^*   \bar P^n \Lambda)\Lambda^*\bar P\Lambda\\
&=&( \bar P^n \Lambda)\Lambda^*\bar P\Lambda\\
&=& \bar P^n (\Lambda\Lambda^*\bar P\Lambda)\\
&=&\bar P^{n+1} \Lambda \eq
The case $n=0$ corresponds to assumption \eqref{H3}.
So we get that for any $\bar A\in\cA(\bar P)$, 
 \bq  \Lambda\Lambda^*\bar A\Lambda&=&\bar A\Lambda\eq
 from which we deduce \eqref{hypK} if $\bar P$ is uniplicit.\wwtbp
\begin{rem}\label{commut}
(a) The inclusion $\cA(P)\cap\cK(V)\subset \cK(P)$ is always true, but it is not necessarily an equality.
Indeed, if $\cK(P)=\cA(P)\cap\cK(V)$, then the elements of $\cK(P)$ commute.
But $\cS_P$ is naturally included into $\cK(P)$ via the representation $\cS_P\ni g\mapsto  T_g\in \cK(V)$
where $T_g$ is the deterministic Markov kernel given by
\bq
\fo x\in V,\qquad T_g(x,\cdot)&=&\delta_{g(x)}\eq
If the elements of $\cK(P)$ commute, then $\cS_P$ is itself commutative.
This is not always true, one can e.g.\ consider the transition kernel of the random walk generated by the transpositions on the permutation group $\cS_N$, with $N\geq 3$.
\par\sm
(b) The example of Remark \ref{tensorization} is equally such that $\cS_{\bar P}$ is not commutative for $N\geq 3$.
Indeed, consider for $\sigma\in\cS_N$ the mapping  $g$ on $\{0,1\}^N$ obtained by shuffling the coordinates according to $\sigma$.
Then $T_g$, defined as above, belongs to $\cS_{\bar P}$. It follows that $\cS_{\bar P}$ contains $\cS_N$ as a subgroup and thus cannot be commutative.
Despite the fact that $\bar P$ is not uniplicit, it was proven in \cite{miclo:hal-01117051} that the conclusion of Proposition \ref{pro4} is true, where $G\df S_N$.
In this case $P$ is a birth and death chain and is thus uniplicit.\par\sm
(c) Even if it outside the finite framework, the example of the Laplacian $L$ on the sphere $\SS^N\subset \RR^{N+1}$, with $N\geq1$, is also such that   $\cK(L)$ (rigorously,
one should define it with respect  to the associated heat kernel at a positive time) is not commutative,
because $\cS_L$ contains all the isometric transformations of $\SS^N$, namely the orthogonal group O($N+1$). 
Note nevertheless that since $\cK(L)$ is big, the same is true for $\cH(L)$: it  is the whole sphere!
We mention this case, because it plays an important role 
 in Carlen, Geronimo and Loss \cite{MR2764893}. At first view, it has some similarities with the situation of (b) above: $L$ is not uniplicit
 but formally the conclusion of Proposition \ref{pro4} is true when  $G$ is the subset of O($N+1$) conserving the norm of the $n$ first coordinates of $\RR^N$,
 with $n\in \lin N-1\rin$.\par\sm
 (d) Despite what we just said, it seems there is an important difference between the cases (b) and (c) above. In the latter
 it can be checked that $\cK(L,\Lambda)\not=\cK(L)$, while in the former we think that $\cK(\bar P,\Lambda)=\cK(\bar P)$.
 That is why Proposition \ref{trop} could be applied to such $\bar P$ without the assumption of uniplicity,
 thus explaining the validity of Proposition \ref{pro4} for this example. In \cite{miclo:hal-01117051}, it was rather used that
 $\cH(\bar P,\Lambda)\not=\emptyset$, where $\cH(\bar P,\Lambda)\df\{x\in\{0,1\}^N\st \delta_x\cK(\bar P,\Lambda)=\cP(\{0,1\}^N)\}$.
\end{rem}
\par
With these observations, we can come to the
\prooff{Proof of Proposition \ref{trop}}
Consider  $x_0\in\cH(\bar P)$. Taking into account \eqref{hypK}, we have
\bq
\delta_{x_0}\Lambda\Lambda^*\cK(\bar P)\Lambda&=&
\delta_{x_0}\cK(\bar P)\Lambda\\
&=&\cP(\bar V)\Lambda\\
&=&\cP(V)\eq
where we used \eqref{H1}.
Assumption \eqref{H2} ensures that $\delta_{x_0}\Lambda=\delta_{\pi(x_0)}$,
so we get
\bq
\delta_{\pi(x_0)}\Lambda^*\cK(\bar P)\Lambda&=&\cP(V)\eq
Finally we use Lemma \ref{CGL} to see that 
\bq \cP(V)&\subset & \delta_{\pi(x_0)}\cK(P)\eq
which is the wanted result.\wwtbp
\par
It is time now to consider
the purely determinist case for $\Lambda$, which simplifies most of the previous hypotheses.
More precisely, assume that there exists a surjective mapping $\pi$ from $\bar V$ to $V$
such that $\Lambda$ is given by
\bqn{deter}
\fo x\in\bar V,\qquad \Lambda(x,\cdot)&\df& \delta_{\pi(x)}(\cdot)\eqn
\begin{lem}\label{hyphyp}
Under \eqref{deter}, if $\bar P$ is a Markov kernel on $\bar V$ such that $(\bar P,\Lambda)$ is positive and if $P$ is a Markov kernel on  $V$ 
satisfying the intertwining relation \eqref{intertw} (called Dynkin's condition in this situation, see \cite{MR0193671}),
then \eqref{hyp}, \eqref{H1}, \eqref{H2} and \eqref{H3} are true. Furthermore, $\Lambda^*\Lambda=I$ and $P$ is given by \eqref{hyp2}.
\end{lem}
\proof
Under Assumption \eqref{deter}, it was seen in Lemma 5 of \cite{miclo:hal-01117051}
that $\Lambda\Lambda^*$ is the conditional expectation with respect to the sigma-algebra $\cT$ generated by $\pi$.
\\
Consider \eqref{hyp}, which amounts to 
\bq
\fo f\in\LL^2(\mu),\qquad 
\Lambda\Lambda^*\bar P\Lambda[f]&=&\bar P\Lambda[f]\eq
Note that the relation $\bar P\Lambda[f]=\Lambda P[f]=P[f]\circ \pi$ implies that $\bar P\Lambda[f]$ is $\cT$-measurable for any $ f\in\LL^2(\mu)$,
so the above equality holds.
Similarly, using that $\Lambda[f]$ is $\cT$-measurable for any $ f\in\LL^2(\mu)$, we get \eqref{H3}.
It follows that $\Lambda^*\Lambda$ is a projection in $\LL^2(\mu)$ and to see that $\Lambda^*\Lambda=I$,
it is sufficient to check that $\Lambda^*\Lambda$ is injective.
So let $f\in\LL^2(\mu)$ be such that $\Lambda^*\Lambda[f]=0$,
we get that \bq
f\circ\pi\ =\ \Lambda[f]\ =\ \Lambda\Lambda^*\Lambda[f]\ =\ 0\eq
Since $\pi$ is surjective, it appears that $f=0$.
\\
It follows that $P$ is given by  \eqref{hyp2}:
\bq
P\ =\ \Lambda^*\Lambda P\ =\  \Lambda^*\bar P\Lambda \eq
Condition \eqref{deter} implies obviously \eqref{H2}, and \eqref{H1} due to the surjectivity of $\pi$.\wwtbp
\par
\begin{rem}\label{nondet}  (a) The deterministic case \eqref{deter} is not the only one where  \eqref{hyp} is satisfied. 
Indeed, assume that $\pi$ is surjective but not injective in \eqref{hyp}. Let  $\bar P$ be a Markov kernel on $\bar V$ such that $(\bar P,\Lambda)$ is positive.
From Lemme \ref{hyphyp}, it appears that $\Lambda\Lambda^*=I$, so we get
\bq
\Lambda\Lambda^*\bar P^*\Lambda^*&=&\bar P^*\Lambda^*\eq
namely  \eqref{hyp} for $\bar P^*$ and $\Lambda^*$.
But since $\pi$ is not injective, the conditional expectation $\Lambda^*\Lambda$ is not the identity, thus $\Lambda^*$ does not satisfy 
\eqref{deter}.\par\sm
(b) Under Assumption \eqref{deter}, Condition \eqref{m0hyp} is also simple to understand: it asks that the conditional expectations with respect to $\cT$ (the sigma-algebra generated by $\pi$)
with respect to $\bar\mu$ and $\bar m_0$ coincide. Namely, if $(A_1, ..., A_l)$ is the partition of $\bar V$ generating $\cT$ (corresponding to the equivalence relation between $x,y\in \bar V$
given by $\pi(x)=\pi(y)$), then
$\bar m_0$ satisfies \eqref{m0hyp} if and only if it is of the form
\bq
\fo x\in\bar V,\qquad m_0(x)&=&\sum_{k\in\lin l\rin}\frac{a_k}{\bar\mu(A_k)}\un_{A_k}(x)\bar\mu(x)\eq
where $(a_1,..., a_l)$ is a probability measure on $\lin l\rin$.
\end{rem}\par
From Proposition \ref{trop} and Lemma \ref{hyphyp}, we deduce:
\begin{cor}\label{cor1}
Assume that the Markov kernel $\bar P$ is uniplicit and that $\cH(\bar P)\not=\emptyset$. 
Let $P$ be a Markov kernel satisfying Relation \eqref{intertw} with a link $\Lambda$ given by \eqref{deter} with $\pi$ surjective.
Then  $P$ is uniplicit and satisfies \eqref{ghp} as well as  the hypergroup property.
\end{cor}
\proof
The above results show that $\cH(P)\not=\emptyset$. According to Lemma \ref{cHgroup}, it is then sufficient to check that $P$ is uniplicit.
By duality, we have $P^*\Lambda^*=\Lambda^*\bar P^*=\Lambda^*\bar P$, it implies, via the equality $\Lambda^*\Lambda=I$ of Lemma~\ref{hyphyp},
\bq
P^*&=&
P^*\Lambda^*\Lambda\\&=&\Lambda^*\bar P\Lambda
\\
&=& P
\eq
where we used \eqref{hyp}, which is true due to Lemma ~\ref{hyphyp} again.
This shows that $P$ is reversible. 
\\
Consider $\theta$ an eigenvector of $P$ and $\varphi ,\wi\varphi $ two associated eigenvectors.
From the intertwining relation \eqref{intertw} we get
\bq
\theta \Lambda[\varphi ]&=&\bar P[\Lambda[\varphi ]]\eq
and similarly for $\wi\varphi $. By uniplicity of $\bar P$, $\Lambda[\varphi ]$ and $\Lambda[\wi\varphi ]$ are then  co-linear.
Remembering that $\Lambda$ is injective by surjectivity of $\pi$, we get that $\varphi $ and $\wi\varphi $ are co-linear
as wanted.\wwtbp
\par
Proposition \ref{pro4} is itself a consequence of the previous corollary. Indeed, it is immediate to check that $\bar P$, $P$ and $\pi$ given before Proposition \ref{pro4} 
satisfy the intertwining relation \eqref{intertw} where $\Lambda$ is defined by \eqref{deter}.\par\me
To end this section, we mention some (upper) semi-continuity properties associated to the Markov commutator convex semi-groups, suggesting the easy handling of this notion.
Note that for any  Markov kernel $P$ on the finite set $V$ and $x\in V$, the sets $\cK(P)$ and $\cK(P,x)$ are compact subsets, respectively of the set of Markov kernels
and of probability measures on $V$ (endowed with the topologies inherited from those of $\RR^{V^2}$ and $\RR^V$), themselves being compact.
As usual, consider the Hausdorff topology on the compact subsets of a compact set, it turns it into a compact set itself.
The following properties are elementary and their proofs are left to the reader.
\begin{lem}\label{uppercont}
Let $(P_n)_{n\in\NN}$ be a sequence of Markov kernels on $V$ converging to $P$. We have for any $x\in V$,
\bq
\limsup_{n\ri\iy} \cK(P_n)&\subset &\cK(P)\\
\limsup_{n\ri\iy} \cK(P_n,x)&\subset &\cK(P,x)\\
\limsup_{n\ri\iy} \cH(P_n)&\subset &\cH(P)\eq
\end{lem}
\par
As a consequence, the set of Markov kernels $P$ on $V$ satisfying the generalized hypergroup property \eqref{ghp} is closed.\par
Let us remark that the above last inclusion can be strict.
Anticipating a little on the next section, consider $V\df\{0,1\}$ and let $(U_n)_{n\in\NN}$ be a sequence of functions on $V$
satisfying $U_n(0)>U_n(1)$ for all $n\in\NN$ and $\lim_{n\ri\iy}U_n=0$.
With the notation of \eqref{MU}, we have
\bq
\fo n\in\NN,\qquad \cH(M_{U_n})&=&\{0\}\\
\lim_{n\ri\iy} M_{U_n}&=&M_0\\
\cH(M_0)&=&\{0,1\}\eq

\section{On the discrete Achour-Trimèche's theorem}

Here the specific birth and death situation is considered in a more detailed way.
The diffusive Achour-Trimèche's theorem will be partially translated into the discrete case,
but first we show it cannot be extended to all convex potentials. It corresponds respectively to the proofs
of Proposition \ref{pro3} and Theorem \ref{theo1}.
\par\me
The previous section provided  all the ingredients necessary to the
\prooff{Proof of Theorem \ref{theo1}}
Recall the setting described in the introduction.
Theorem \ref{theo1} is proven by a contradictory argument: assume there exists a generalized Metropolis procedure
$\cC\ni U\mapsto P_U$ such that $P_U$ satisfies the hypergroup property for all $U\in\cC$.\\
Since $N\geq 2$, there exists $U\in\cC$ such that $U(0)=U(1)$ and which is not symmetric with respect to the mapping
$\lin 0, N\rin \ni x\mapsto N-x$.
For $\epsilon>0$, consider the function $U_{\epsilon}$ defined on $\lin 0,N\rin$ by
\bq
\fo x\in\lin 0,N\rin,\qquad 
U_\epsilon(x)&\df&
\lt\{
\begin{array}{ll}
U(0)+\epsilon&\hbox{, if $x=0$}\\
U(x)&\hbox{, otherwise}\end{array}\rt.\eq
It is clear that $U_\epsilon\in\cC$. 
Furthermore, due to the convexity of $U$ and the assumption $U(0)=U(N)$,
it appears that 
$U_\epsilon(0)>U_\epsilon(x)$ for all $x\in \lin N\rin$.
By Definition \eqref{muU}, the minimum of $\mu_{U_\epsilon}$ is only attained at 0.
Taking into account Proposition \ref{pro1}, it follows that $\cH(P_{U_\epsilon})=\{0\}$.
By letting $\epsilon>0$ go to zero,
Lemma \ref{uppercont} implies that 
$0\in\cH(P_U)$.
The same reasoning, where the value of $U(N)$ is a little increased,  equally enables to conclude that $N\in\cH(P_U)$.
So we get that $\{0,N\}\subset \cH(P_U)$.
Since $P_U$ is a birth and death, it is uniplicit, and according to Proposition \ref{pro2}, we can find $g\in\cS_{P_U}$
with $g(0)=N$.
Note that under the action of any element of the symmetry group $\cS_P$, the graph of the transitions permitted  by $P$
is preserved (not taking into account the self-loops). For  birth and death transitions on $\lin 0,N\rin$, this graph is the usual linear graph structure of $\lin 0,N\rin$.
There are only two graph morphisms preserving this structure, the identity and the mapping $\lin 0, N\rin \ni x\mapsto N-x$.
So we end up with a contradiction, because $g$ can be neither of them.\wwtbp
\par\me
We now come to the proof of Proposition \ref{pro3}. We begin by reducing the problem to symmetric potentials.
Recall that the classical Metropolis procedure $\cC\ni U\mapsto M_U$ is defined by \eqref{MU}.
\begin{lem}\label{desymetrisation}
If for all $N\in \NN\setminus\{1\}$, the Metropolis kernel $M_U$ satisfies the hypergroup property for $U\in\wi\cC_{\mathrm{s}}$,
then it is also true for $U\in \wi\cC_{\mathrm{m}}$.
\end{lem}
\proof
This is a consequence of Proposition
\ref{pro4}.
Indeed, let $U\in\wi\cC_{\mathrm{m}}$, up to reversing the discrete segment $\lin 0,N\rin$, assume that $U$ is non-increasing. Consider $\bar V\df \lin 0,2N+1\rin$, on which we construct the potential $\bar U$ by symmetrization of $U$ with respect to $N+1/2$.
 Note that $\bar U$ is convex
 and more precisely that $\bar U\in \wi\cC_{\mathrm{m}}$, due to the assumption $U(N-1)-U(N)\geq 2\ln(2)$, which implies 
 \bq
\bar U(N+2)-\bar U(N+1)&\geq& 2\ln(2)\\
&=&\bar U(N+1)-\bar U(N)+2\ln(2)\\
&\geq & \bar U(N)-\bar U(N-1)+4\ln(2)\eq
Associate to $\bar U$ the classical Metropolis kernel $\bar M_{\bar U}$ on $\bar V$.
Let $\bar G=\cS_{\bar M_{\bar U}}$ be the group consisting of the identity and of the involution $\lin 0, 2N+1\rin \ni x\mapsto 2N+1-x$.
The reduction presented before Proposition
\ref{pro4} transforms $\bar M_{\bar U}$ into $M_U$ (up to a modification of the constant $\Sigma_U$ given in \eqref{SigmaU}, which has no impact on the hypergroup property,
since it amounts to change $M_U$ into a convex combination of $M_U$ and $I$). Again, since $\bar M_{\bar U}$ is a birth and death chain, it is uniplicit. Thus Proposition
\ref{pro4} enables to see that $M_U$ satisfies Condition \eqref{ghp}, because by assumption this is true for $\bar M_{\bar U}$.
Applying once more Lemma \ref{cHgroup}  shows that $M_U$ satisfies the hypergroup property.\wwtbp
\begin{rem}\label{wiM0}
In the above proof, another symmetrization could have been considered: let $\bar V\df  \lin 0,2N\rin$ and $\bar U$ be obtained from $U$ by symmetry with respect to $N$ ($U$ being non-increasing). Applying the same arguments under the relaxed assumption 
$U(N-1)-U(N)\geq \ln(2)$ (implying $\bar U(N+1)-\bar U(N)\geq \bar U(N)-\bar U(N-1)+2\ln(2)$) , we get in the end that $\wi M_U$ satisfies the hypergroup property,
where $\wi M_U$ is defined as $M_U$ in \eqref{MU}, but with $M_0$ replaced by the exploration kernel $\wi M_0$ given by 
\bq
\fo x\not=y\in\lin 0,N\rin,\qquad\wi M_0(x,y)&\df&\lt\{
\begin{array}{ll}
1/2&\hbox{, if $\lve x-y\rve=1$ and $x\not=N$}\\
1&\hbox{, if $x=N$ and $y=N-1$}\\
0&\hbox{, otherwise}
\end{array}\rt.
\eq
\end{rem}
\par
It remains to prove that for $U\in\wi\cC_{\mathrm{s}}$, $M_U$ satisfies the hypergroup property. We did not find general arguments to obtain this result.
Instead, we will adapt to the discrete case the proof presented by Bakry and Huet \cite{MR2483738} in the context of symmetric one-dimensional diffusions.
\begin{pro}\label{cCs}
For any $U\in\wi\cC_{\mathrm{s}}$, 
the Metropolis kernel $M_U$ satisfies the hypergroup property, with respect to the points 0 and $N$.
\end{pro}
By uniplicity of $M_U$ and its symmetry with respect to the mapping 
\bqn{s} s\st \lin 0, N\rin\ni x&\mapsto& N-x\eqn
it is sufficient to check that $0\in\cH(M_U)$ for given $U\in\wi\cC_{\mathrm{s}}$.
Let us consider more generally the problem of showing that $0\in\cH(P)$, when $P$ is an irreducible birth and death Markov transition
on $\lin 0,N\rin$, left invariant by the symmetry $s$.
By definition, it amounts to show that for any given probability $m_0\in\cP(\lin 0,N\rin)$, there is a Markov kernel $K$ commutating with $P$ and such that
$K(0,\cdot)=m_0$.
This question is equivalent to the fact that a wave equation starting from a non-negative condition remains non-negative,
as it was shown by Bakry and Huet \cite{MR2483738} in the diffusive situation and in Remark 6 of \cite{miclo:hal-01117051} for the discrete case.
More precisely, there is a unique matrix $K$ commuting with $P$ such that $K(0,\cdot)=m_0$ (due to the uniplicity of $M_U$, see the proof of Lemma \ref{unic} or Lemma 10 of 
\cite{miclo:hal-01117051}), our problem is to check that its entries are non-negative.
Denote $L=P-I$, the Markovian generator matrix associated to $P$ and
\bq
 \fo x,y\in\lin 0, N\rin,\qquad
 k(x,y)&\df& \frac{K(x,y)}{\mu(y)}\eq
 The commutation of $K$ with $P$ can be rewritten as
 the wave equation
 \bqn{wave}
  \fo x,y\in \lin 0, N\rin,\qquad
  L^{(1)}[k](x,y)&=&  L^{(2)}[k](x,y)
\eqn
where for $i\in\{1,2\}$, $L^{(i)}$ stands for the generator acting on the $i$-th variable as $L$.\\
Consider the discrete triangle
\bq
\tr&\df& \{(x,y)\in \lin 0,N\rin^2\st x\leq y\hbox{ and } x\leq N-y\}\eq
For $z_0\df(x_0,y_0)\in \tr$, let $p_{z_0}^{-}\df(p_{z_0}^{-}(n))_{n\in \lin 0, 2y_0\rin}$ be the path defined by iteration through
\bq
p_{z_0}^{-}(0)&\df&z_0\\
\fo n\in  \lin 0, 2y_0-1\rin,\qquad
p_{z_0}^{-}(n+1)&\df&
\lt\{
\begin{array}{ll}
p_{z_0}^{-}(n)-(0,1)&\hbox{, if $n$ is even}\\
p_{z_0}^{-}(n)-(1,0)&\hbox{, if $n$ is odd}\end{array}\rt.
\eq
Note that the path $p_{z_0}^{-}$ stays in $\tr$ and that $p_{z_0}^{-}(2y_0)$ belongs to the segment $\lin 0,N\rin\times\{0\}$.
\\
Similarly, for $z_0\in \tr$, we define the path $p_{z_0}^{+}\df(p_{z_0}^{+}(n))_{n\in \lin 0, 2y_0\rin}$, which is symmetric to $p_{z_0}^{-}$
with respect to the axe $x=x_0$.
The interest of these paths is:
\begin{lem}\label{ipp3}
Assume that the mapping $k\st \lin 0, N\rin^2\ri \RR$ satisfies the wave equation \eqref{wave}.
Then for any $z_0\df(x_0,y_0)\in \tr$, we have, if $y_0\geq 1$,
\bq
\omega(z_0,p_{z_0}^-(1))k(z_0)&=&
[\omega(z_0,p_{z_0}^-(1))-\omega(p_{z_0}^-(1),p_{z_0}^-(2))-\omega(p_{z_0}^+(1),p_{z_0}^+(2))]k(p_{z_0}^-(1))\\
&&+
\omega(p_{z_0}^{-}(2y_0-1),p_{z_0}^{-}(2y_0))k(p_{z_0}^{-}(2y_0))+
\omega(p_{z_0}^{+}(2y_0-1),p_{z_0}^{+}(2y_0))k(p_{z_0}^{+}(2y_0))
\\&&
+\sum_{n\in \lin 2, 2y_0-1\rin}[\omega(p_{z_0}^{-}(n-1),p_{z_0}^{-}(n))-\omega(p_{z_0}^{-}(n),p_{z_0}^{-}(n+1))]k(p_{z_0}^{-}(n))\\&&
+\sum_{n\in \lin 2, 2y_0-1\rin}[\omega(p_{z_0}^{+}(n-1),p_{z_0}^{+}(n))-\omega(p_{z_0}^{+}(n),p_{z_0}^{+}(n+1))]k(p_{z_0}^{+}(n))
\eq
where for any $(z,z')\df((x,y),(x',y'))\in \lin 0,N\rin^4$, we take
\bq
\omega(z,z')&\df&\lt\{
\begin{array}{ll} \mu(x)\mu(y)L(x,x')&\hbox{, if $z'-z\in \{(1,0), (-1,0)\}$}\\
\mu(x)\mu(y)L(y,y')&\hbox{, if $z'-z\in \{(0,1), (0,-1)\}$}
\end{array}\rt.
\eq
\end{lem}
\proof
From the reversibility of $L$ with respect to $\mu$, we deduce the discrete integration by part formula: for any functions $f,g$ on $\lin 0, N\rin$,
we have
\bq
\mu[fL[g]]&=&-\sum_{0\leq x<y\leq N} \mu(x)L(x,y)[f(y)-f(x)][g(y)-g(x)]\eq
In particular, if $f$ is the indicator function of a segment $\lin q,r\rin\subset \lin 0,N\rin$, we get
\bqn{ipp}
\mu[\un_{\lin q,r\rin}L[g]]&=& [g(r+1)-g(r)]\mu(r)L(r,r+1)+ [g(q-1)-g(q)]\mu(r)L(q,q-1)\eqn
with the convention (Neumann boundary) that $g(-1)=g(0)$ and $g(N+1)=g(N)$.\\
For $z_0\in \tr$, define the discrete triangle
\bqn{trz0}
\tr(z_0)&\df& \{(x,y)\in \lin 0,N\rin^2\st x\leq y-y_0+x_0-1\hbox{ and } x\leq -y+y_0+x_0-1\}\eqn
Applying \eqref{ipp} horizontally and vertically, we get, for $k$  satisfying the wave equation \eqref{wave},
\bqn{ipp2}
\nonumber 0&=&\mu^{\otimes 2}[\un_{\tr(z_0)}(L^{(1)}-L^{(2)})[k]]\\
&=&\sum_{e\in\pa \tr(z_0)} dk(e) \chi(e) \omega(e)\eqn
where the boundary $\pa\tr(z_0)$ of $\tr(z_0)$ is defined by
\bq
\pa\tr(z_0)&\df& \{(z,z')\in \tr(z_0)\times (\lin 0,N\rin^2\setminus \tr(z_0))\st z'-z\in\{(1,0), (-1,0), (0,1), (0,-1)\}\}\eq
and where for any $e\df(z,z')\in \pa\tr(z_0)$, $\omega(e)$ was defined in the statement of the lemma and
\bq
dk(e)&\df& k(z')-k(z)\\
\chi(z,z')&\df&\lt\{
\begin{array}{ll} 1&\hbox{, if $z'-z\in \{(1,0), (-1,0)\}$}\\
-1&\hbox{, if $z'-z\in \{(0,1), (0,-1)\}$}
\end{array}\rt.
\eq
It is easy (but a picture can help) that \eqref{ipp2}
can written under the form 
\bqn{bord}
\nonumber 0
&=&\sum_{n\in \lin 0, 2y_0-1\rin} [k(p_{z_0}^{-}(n+1))-k(p_{z_0}^{-}(n))]\omega(p_{z_0}^{-}(n),p_{z_0}^{-}(n+1))\\
&&+\sum_{n\in \lin 1, 2y_0-1\rin}[ k(p_{z_0}^{+}(n+1))-k(p_{z_0}^{+}(n))]\omega(p_{z_0}^{+}(n),p_{z_0}^{+}(n+1))\eqn
Observe that the first sum can be transformed (via discrete integration by parts, also known as Abel's trick) into
\bq
\lefteqn{\sum_{n\in \lin 0, 2y_0-1\rin} [k(p_{z_0}^{-}(n+1))-k(p_{z_0}^{-}(n))]\omega(p_{z_0}^{-}(n),p_{z_0}^{-}(n+1))}\\
&=&
k(p_{z_0}^{-}(2y_0))\omega(p_{z_0}^{-}(2y_0-1),p_{z_0}^{-}(2y_0))-k(z_0)\omega(z_0,p_{z_0}^{-}(1))\\&&
-\sum_{n\in \lin 1, 2y_0-1\rin}k(p_{z_0}^{-}(n))[\omega(p_{z_0}^{-}(n),p_{z_0}^{-}(n+1))-\omega(p_{z_0}^{-}(n-1),p_{z_0}^{-}(n))]\eq
A similar manipulation is possible for the second sum \eqref{bord} and we end up with the result announced in the lemma.\wwtbp
As a consequence, we get
\begin{pro}\label{trtr}
Assume that $P$ is a birth and death transition kernel on $\lin 0,N\rin$  such that 
\bq
\fo z\df(x,y)\in\tr,\qquad
P(y-1,y)&\geq & P(x,x-1) +P(x,x+1)\\
\fo z\df(x,y)\in\wi \tr,\qquad
P(y,y-1)&\leq & P(x-1,x) \wedge P(x+1,x)
\eq
where $\wi\tr$ is the ``interior" of $\tr$:
\bq
\wi \tr&\df&  \{(x,y)\in \lin 0,N\rin^2\st x\leq y-1\hbox{ and } x\leq N-y-1\}\eq
Let $k$ be a solution of \eqref{wave} such that $k(\cdot,0)$ is non-negative. Then $k$ remains non-negative on $\tr$.
\end{pro}
\proof
We begin by showing that the condition of the proposition (which can be written identically in terms of $L$), 
implies that
for any $z_0\df(x_0,y_0)\in\tr$ and $n\in  \lin 2, 2y_0-1\rin$, we have
\bq
\omega(p_{z_0}^{-}(n-1),p_{z_0}^{-}(n))-\omega(p_{z_0}^{-}(n),p_{z_0}^{-}(n+1))&\geq &0\\
\omega(p_{z_0}^{+}(n-1),p_{z_0}^{+}(n))-\omega(p_{z_0}^{+}(n),p_{z_0}^{+}(n+1))&\geq &0\eq
It amounts to see that for any $(x,y)\in\tr$,
\bqn{ww1}
\lt\{
\begin{array}{rcl}\omega((x,y),(x,y-1))-\omega((x,y-1),(x-1, y-1))&\geq &0\\
\omega((x,y),(x,y-1))-\omega((x,y-1),(x+1, y-1))&\geq &0
\end{array}\rt.
\eqn
and that for any $(x,y)\in\wi\tr$,
\bqn{ww2}
\lt\{
\begin{array}{rcl}
\omega((x,y),(x-1,y))-\omega((x-1,y),(x-1, y-1))&\geq &0\\
\omega((x,y),(x+1,y))-\omega((x+1,y),(x+1, y-1))&\geq &0
\end{array}\rt.
\eqn
Concerning \eqref{ww1}, let $\varepsilon\in\{-1,+1\}$,
we have
\bq
\lefteqn{\omega((x,y),(x,y-1))-\omega((x,y-1),(x+\varepsilon, y-1))}\\&=&\mu(x)\mu(y)L(y,y-1)-\mu(x)\mu(y-1)L(x,x+\varepsilon)\\
&=&\mu(x)\mu(y-1)[L(y-1,y)-L(x,x+\varepsilon)]
\eq
where we used the reversibility of $\mu$ with respect to $L$.
By the first assumed inequality, we have in particular
$P(y-1,y)\geq P(x,x-1) \vee P(x,x+1)$, so that
the last r.h.s.\ is non negative, as wanted.
\\
The treatment of \eqref{ww2} is similar, taking into account the second assumed inequality:
\bq
\lefteqn{\omega((x,y),(x+\varepsilon,y))-\omega((x+\varepsilon,y),(x+\varepsilon, y-1))}\\&=&\mu(x)\mu(y)L(x,x+\varepsilon)-\mu(x+\varepsilon)\mu(y)L(y,y-1)\\
&=&\mu(x+\varepsilon)\mu(y)[L(x+\varepsilon,x)-L(y,y-1)]\\
&\geq &0
\eq
Next we want to show that
\bq\omega(z_0,p_{z_0}^-(1))-\omega(p_{z_0}^-(1),p_{z_0}^-(2))-\omega(p_{z_0}^+(1),p_{z_0}^+(2))&\geq &0\eq
Writing $(x,y)\df p_{z_0}^-(1)$, it means that
\bq
\omega((x,y+1),(x,y))-\omega((x,y),(x-1,y))-\omega((x,y),(x+1,y))&\geq &0\eq
namely
\bq
\mu(x)\mu(y)[L(y,y+1)-L(x,x+1)-L(x,x-1)]&\geq &0\eq
condition which is satisfied by the first assumed inequality of the lemma (since $z_0=(x,y+1)$).\\
Thus all the coefficients in front of values of $k$ in the equality of Lemma \ref{ipp3} are non-negative. Assume that $k$ does not remain non-negative on $\tr$.
We can then consider $y_0$ the minimal value of $y\in \lin 0,N\rin$ such that there exists $y\leq x\leq N-x$ such that $k(x,y)<0$.
Next, let $x_0$ the minimal value of $x\in\lin y,N-y\rin$ such that $k(x,y_0)<0$.
In particular, $z_0\df(x_0,y_0)\in \tr$ and $k(z_0)<0$, fact which is in contradiction with the equality of Lemma \ref{ipp3}, whose r.h.s.\ is non-negative.\wwtbp
\par
Assume now that $P$ is furthermore left invariant by the symmetry $s$ defined in \eqref{s}.
One important consequence is that the conclusion of Proposition \ref{trtr} is valid on the whole discrete square $\lin 0, N\rin^2$:
\begin{pro}\label{carre}
Assume that the birth and death transition $P$ on $\lin 0,N\rin$ is invariant by $s$. Let $k$ be a solution of \eqref{wave}.
Then $k$ is left invariant by the following symmetries of the discrete square:
\bq
\lin 0, N\rin^2\ni(x,y)&\mapsto&  (y,x)\\
\lin 0, N\rin^2\ni(x,y)&\mapsto&  (N-x,N-y)\\
\lin 0, N\rin^2\ni(x,y)&\mapsto&  (N-y,N-x)\eq
As a consequence, if $k$ is non-negative on $\tr$, then it is non-negative on $\lin 0, N\rin^2$.
\end{pro}
\proof
Consider $\wi k\st \tr\ri\RR$ satisfying the wave equation \eqref{wave} on $\wi \tr$.
Extend $\wi k$ to the discrete triangle $\tr_2\df\{(x,y)\in \lin 0,N\rin^2\st y\leq N-x\}$ by symmetry with respect to the line $y=x$.
Let us check that $\wi k$ satisfies \eqref{wave} on $\wi \tr_2\df \{(x,y)\in \lin 0,N\rin^2\st y\leq N-x-1\}$.
By symmetry of $P$, it is obvious on the image of $\wi\tr$ by the mapping $(x,y)\mapsto (y,x)$.
Thus it is sufficient to show that
\eqref{wave} is also valid on the points $(x,x)\in \wi\tr_2$. Indeed, we compute that
\bq
\lefteqn{L^{(1)}[\wi k](x,x)-L^{(2)}[\wi k](x,x)}\\
&=& L(x,x+1)(\wi k(x+1,x)-\wi k(x,x))+L(x,x-1)(\wi k(x-1,x)-\wi k(x,x))\\&&
-L(x,x+1)(\wi k(x,x+1)-\wi k(x,x))-L(x,x-1)(\wi k(x,x)-\wi k(x,x-1))\\
&=&0\eq
due to the construction by symmetrization. \\
Next we can extend $\wi k$ to $\lin 0,N\rin^2$ by symmetrization with respect to the line $y=N-x$.
The same arguments as above show that this extension satisfies \eqref{wave} on $\lin 0,N\rin^2$.
Observe that the mapping $\wi k$ constructed in this way is left invariant by the symmetries presented in the lemma.
\par
Now consider $k\st \lin 0,N\rin^2\ri\RR$ a solution of \eqref{wave}. Let $\wi k$ be its restriction to $\tr$.
By the above construction, we extend $\wi k$ to $\lin 0,N\rin^2$ into a function also satisfying \eqref{wave}.
Note that  $k(\cdot,0)=\wi k(\cdot,0)$, so by uniqueness of the solution of \eqref{wave}
given its value on the discrete segment $\{0\}\times \lin 0, N\rin$, we get $k=\wi k$.\wwtbp\par
Consider the following assumption called (H): 
 the mappings $\lin 0, \lfloor N/2\rfloor \rin\ni x\mapsto 2^xP(x,x+1)$ and $\lin 0,   \lfloor N/2\rfloor\rin\ni x\mapsto P(x+1,x)$ are respectively  non-increasing and
 non-decreasing.\par
 Our main result about a partial extension of Achour-Trimèche's theorem to the discrete setting can be stated as
\begin{theo}\label{critere}
 Assume that the birth and death transition $P$ on $\lin 0,N\rin$ is invariant by $s$ and that (H) is fulfilled.
 Then $P$ satisfies the hypergoup property with respect to $0$ and $N$.
 \end{theo}
 \proof
 According to Proposition \ref{carre}, it is enough to check that (H) implies the assumption of Proposition~\ref{trtr}.
 Note that in the case where $N$ is odd, by symmetry of $P$ through $s$, we have
 $P((N-1)/2,(N+1)/2)=P((N+1)/2,(N-1)/2)$. When $N$ is even, we rather get $P(N/2,N/2+1)=P(N/2,N/2-1)$ and
 $P(N/2-1,N/2)=P(N/2+1,N/2)$.
In both situations, it appears that (H) leads to \bq
\fo y\in\lin 0, \lfloor N/2\rfloor-1\rin,\,\fo x\in \lin y+1, \lfloor N/2\rfloor\rin,\qquad
\lt\{\begin{array}{c}2P(x+1,x)\leq 2P(x,x+1)\leq P(y,y+1)\\
P(y+1,y)\leq P(x+1,x)\leq P(x,x+1)
\end{array}\rt.\eq
By symmetry of $P$ through $s$, it follows that
\bq
\fo y\in\lin 0, \lfloor N/2\rfloor-1\rin,\,\fo x\in \lin y+1, N-y-1\rin,\qquad
\lt\{\begin{array}{c}
 P(x+1,x)+P(x,x+1)\leq P(y,y+1)\\
P(y+1,y)\leq P(x+1,x)\wedge  P(x,x+1)
\end{array}\rt.
\eq
which is the assumption of Proposition \ref{trtr}.
\wwtbp
\par
As a simple corollary we obtain Proposition \ref{cCs}, because $P_U$ satisfies (H) if $U\in\wi\cC_{\mathrm{s}}$.
Indeed, this condition asks for
 the mappings $\lin 0, \lfloor N/2\rfloor \rin\ni x\mapsto U(x)-U(x+1)+2\ln(2)x$ and $\lin 0,   \lfloor N/2\rfloor\rin\ni x\mapsto U(x+1)-U(x)$ to be respectively  non-increasing and
 non-decreasing. This is valid, by the definition of  $\wi \cC$ given before Proposition \ref{pro3}.\par
\begin{rems}\label{finale}
(a) One can replace the exploration kernel $M_0$ given in \eqref{M0} by $\wit M_0$ defined via
\bq
\fo x\not=y\in\lin 0,N\rin,\qquad\wit M_0(x,y)&\df&\lt\{
\begin{array}{ll}
1/2&\hbox{, if $\lve x-y\rve=1$, $x\not=0$ and $x\not=N$}\\
1&\hbox{, if $(x,y)=(0,1)$ or $(x,y)=(N,N-1)$}\\
0&\hbox{, otherwise}
\end{array}\rt.
\eq
The corresponding Metropolis procedure $\wi\cC_{\mathrm{s}}\ni U\mapsto \wit M_U$
(where $\wit M_U$ is defined as in \eqref{MU}, with $M_0$ replaced by $\wit M_0$) also satisfies the hypergroup property,
because (H) is equally true for these birth and death Markovian transitions.
\\
Taking into account Remark \ref{wiM0}, this result can be extended to the Metropolis procedure $\wi\cC_{\mathrm{m}}\cup \wi\cC_{\mathrm{s}}\ni U\mapsto \wit M_U$.\par
Nevertheless, due to the fact that $0\not\in \wi \cC$, we are not able to recover that $\wit M_0$ satisfies the hypergroup property, as it was shown  in Example 7 of \cite{miclo:hal-01117051}.
\par\sm
(b) For $U\in\cC$, consider the  variant classical Metropolis procedure $\wideparen{M}_U$ given by
\bqn{wMU}
\fo x\not=y\in\lin 0,N\rin,\qquad \wideparen{M}_U(x,y)&\df& \wit M_0(x,y)\exp(-(U(y)-U(x))_+)\eqn
Simulations suggest that $\wideparen{M}_U$ satisfies the hypergroup property if the convex function $U$ is either
monotonous or symmetric with respect to the middle point of the discrete segment $\lin 0,N\rin$.
It would be a nice discrete extension of the Achour-Trimèche's theorem, but we have not been able to prove this conjecture.\par\sm
(c) The previous conjecture is not true if in \eqref{wMU}, $\wit M_0$ is replaced by $M_0$ (given by \eqref{M0}).
Indeed, consider the case $N=2$ and $U=0$. Let $k$ be the solution of the corresponding wave equation \eqref{wave}
starting from $k(\cdot,0)\df(0,1,0)$.
Equation  \eqref{wave} at point $(1,1)$ writes:
\bq
\frac12(k(0,0)-k(1,0))+\frac12(k(0,0)-k(1,0))&=&\frac12(k(1,1)-k(1,0))\eq
namely $k(1,1)=-k(1,0)=-1$.
So non-negativity is not preserved by \eqref{wave} and by consequence $M_0$ does not satisfy the hypergroup property.\par
In particular the assumption $U\in\wi\cC$ is not merely technical in Proposition \ref{pro3}. Note this observation is not in contradiction
with the conjecture given in (b).\par\sm
(d) Theorem \ref{critere} enables to construct other examples of birth and death Metropolis procedures satisfying the hypergroup property.
E.g.\ consider the exploration kernel $\widecheck M_0$ given by
\bq
\fo x\not=y\in\lin 0,N\rin,\qquad\widecheck M_0(x,y)&\df&\lt\{
\begin{array}{ll}
1/2^{x\wedge (N-x)}&\hbox{, if $\lve x-y\rve=1$}\\
0&\hbox{, otherwise}
\end{array}\rt.
\eq
Let $\widecheck \cC_{\mathrm{s}}$ be the set of potentials $U$ symmetric with respect to $N/2$ and such that $\widecheck U\in\cC$, where
\bq
\fo x\in\lin 0, N\rin,\qquad \widecheck U(x)&\df& U(x)+\ln(2) (x\wedge (N-x))\eq
Define the Markov kernel $\widecheck{M}_U$ via \bq
\fo x\not=y\in\lin 0,N\rin,\qquad \widecheck{M}_U(x,y)&\df& \widecheck M_0(x,y)\exp(-(\widecheck U(y)-\widecheck U(x))_+)\eq
For $U\in\widecheck \cC_{\mathrm{s}}$, $\widecheck{M}_U$ satisfies (H) and admits $\mu_U$, the Gibbs measure defined in \eqref{muU}, as reversible measure. Thus
$ \cC_{\mathrm{s}}\ni U\mapsto \widecheck{M}_U$ is a generalized birth and death Metropolis procedure satisfying the hypergroup property.
The proof of Lemma \ref{desymetrisation} enables to deduce a similar construction for monotonous potentials
(for instance for convex potentials $U$ such that $\lin 0,N\rin \ni x\mapsto U(x)+\ln(2)x$ is  non-increasing).\\
Note that the potentials from $\widecheck \cC_{\mathrm{s}}$ are more general than those from $\wit \cC_{\mathrm{s}}$,
since the former ones can grow linearly (away from the middle point of the state space), while the latter ones must grow quadratically. The drawback is that  $\widecheck{M}_U$
is further away  from the continuous model $\pa^2-U'\pa$ than $M_U$ defined in \eqref{MU}.
\end{rems}
\par
To finish, let us mention a non-negativity preservation on edges rather than on vertices under a
natural relaxation of the assumption of Proposition \ref{trtr}:
\begin{pro}\label{edges}
Assume that $P$ is a birth and death transition kernel on $\lin 0,N\rin$  such that 
\bq
\fo z\df(x,y)\in\tr,\qquad
P(y-1,y)&\geq & P(x,x-1) \vee P(x,x+1)\\
\fo z\df(x,y)\in\wi \tr,\qquad
P(y,y-1)&\leq & P(x-1,x) \wedge P(x+1,x)
\eq
Let $k$ be a solution of \eqref{wave} such that $k(\cdot,0)$ is non-negative.
Then for any $(x,y)\in\tr$, we have
\bq
(x,y+1)\in\tr&\Rightarrow& k(x,y)+k(x,y+1)\geq 0\\
(x+1,y)\in\tr 
&\Rightarrow& \mu(x)k(x,y)+\mu(x+1)k(x+1,y)\geq 0\eq
\end{pro}
\proof
Note that the equality of Lemma \ref{ipp3} can be rewritten under the form:
\bq
\lefteqn{\omega(z_0,p_{z_0}^-(1))[k(z_0)+k(p_{z_0}^-(1))]}\\&=&
\omega(p_{z_0}^{-}(2y_0-1),p_{z_0}^{-}(2y_0))k(p_{z_0}^{-}(2y_0))+
\omega(p_{z_0}^{+}(2y_0-1),p_{z_0}^{+}(2y_0))k(p_{z_0}^{+}(2y_0))
\\&&
+\sum_{n\in \lin 1, 2y_0-1\rin}[\omega(p_{z_0}^{-}(n-1),p_{z_0}^{-}(n))-\omega(p_{z_0}^{-}(n),p_{z_0}^{-}(n+1))]k(p_{z_0}^{-}(n))\\&&
+\sum_{n\in \lin 1, 2y_0-1\rin}[\omega(p_{z_0}^{+}(n-1),p_{z_0}^{+}(n))-\omega(p_{z_0}^{+}(n),p_{z_0}^{+}(n+1))]k(p_{z_0}^{+}(n))
\eq
So the first implication of the above proposition can be shown as in the proof of Proposition \ref{trtr}, which enables
to see that $k(z_0)+k(p_{z_0}^-(1))\geq 0$, if $y_0\geq 1$.\par
For the second implication, rather consider for $z_0\in\tr$ such that $z_0+(1,0)\in\tr$,
the path $p_{z_0}^{+}\df(p_{z_0}^{+}(n))_{n\in \lin 0, 2y_0+1\rin}$  defined by iteration through
\bq
p_{z_0}^{+}(0)&\df&z_0\\
\fo n\in  \lin 0, 2y_0\rin,\qquad
p_{z_0}^{+}(n+1)&\df&
\lt\{
\begin{array}{ll}
p_{z_0}^{+}(n)+(1,0)&\hbox{, if $n$ is even}\\
p_{z_0}^{+}(n)-(0,1)&\hbox{, if $n$ is odd}\end{array}\rt.
\eq
The set $\tr(z_0)$ defined in \eqref{trz0} must be modified into the ``almost triangle"
\bq
\tr(z_0)&\df& \{(x,y)\in \lin 0,N\rin^2\st x\leq y-y_0+x_0-1\hbox{ and } x\leq -y+y_0+x_0\}\eq
The proof of Lemma  \ref{ipp3} then leads to
\bq
\lefteqn{\omega(z_0,p_{z_0}^-(1))k(z_0)+\omega(p_{z_0}^+(1),p_{z_0}^+(2))k(p_{z_0}^+(1))}\\&=&
\omega(p_{z_0}^{-}(2y_0-1),p_{z_0}^{-}(2y_0))k(p_{z_0}^{-}(2y_0))+
\omega(p_{z_0}^{+}(2y_0),p_{z_0}^{+}(2y_0+1))k(p_{z_0}^{+}(2y_0+1))
\\&&
+\sum_{n\in \lin 1, 2y_0-1\rin}[\omega(p_{z_0}^{-}(n-1),p_{z_0}^{-}(n))-\omega(p_{z_0}^{-}(n),p_{z_0}^{-}(n+1))]k(p_{z_0}^{-}(n))\\&&
+\sum_{n\in \lin 2, 2y_0\rin}[\omega(p_{z_0}^{+}(n-1),p_{z_0}^{+}(n))-\omega(p_{z_0}^{+}(n),p_{z_0}^{+}(n+1))]k(p_{z_0}^{+}(n))
\eq
The proof of Proposition \ref{trtr} now implies that 
$\omega(z_0,p_{z_0}^-(1))k(z_0)+\omega(p_{z_0}^+(1),p_{z_0}^+(2))k(p_{z_0}^+(1))\geq 0$, namely
$\mu(x_0)k(z_0)+\mu(x_0+1)k(x_0+1,y_0)\geq 0$.\wwtbp
\par
The advantage of Proposition \ref{edges} over  Proposition \ref{trtr} is that it enables to recover by approximation (with $N$ going to infinity) the  result of Bakry and Huet \cite{MR2483738}
concerning the preservation of non-negativity by the wave equation in the context of the diffusive Achour-Trimèche theorem.

\bigskip\par\hskip5mm\textbf{\large Acknowledgments:}\par\sm\noindent 
This paper was motivated by the conjecture of Dominique Bakry that the Achour-Trimèche theorem
would be true for all convex potentials in the one-dimensional diffusive case. I'm very grateful to him for all the discussions we had on the subject.
I'm also thankful to  the ANR STAB (Stabilité du comportement asymptotique d'EDP, de processus stochastiques et de leurs discrétisations) for its support.

\vskip2cm
\hskip70mm\box5


\begin{thebibliography}{10}

\bibitem{MR552062}
Abdennebi Achour and Khalifa Trimeche.
\newblock Op\'erateurs de translation g\'en\'eralis\'ee associ\'es \`a un
  op\'erateur diff\'erentiel singulier sur un intervalle born\'e.
\newblock {\em C. R. Acad. Sci. Paris S\'er. A-B}, 288(7):A399--A402, 1979.

\bibitem{MR876954}
David Aldous and Persi Diaconis.
\newblock Strong uniform times and finite random walks.
\newblock {\em Adv. in Appl. Math.}, 8(1):69--97, 1987.

\bibitem{MR2483738}
D.~Bakry and N.~Huet.
\newblock The hypergroup property and representation of {M}arkov kernels.
\newblock In {\em S\'eminaire de probabilit\'es {XLI}}, volume 1934 of {\em
  Lecture Notes in Math.}, pages 295--347. Springer, Berlin, 2008.

\bibitem{MR2764893}
Eric~A. Carlen, Jeffrey~S. Geronimo, and Michael Loss.
\newblock On the {M}arkov sequence problem for {J}acobi polynomials.
\newblock {\em Adv. Math.}, 226(4):3426--3466, 2011.

\bibitem{MR0374545}
Houcine Chebli.
\newblock Op\'erateurs de translation g\'en\'eralis\'ee et semi-groupes de
  convolution.
\newblock In {\em Th\'eorie du potentiel et analyse harmonique ({J}ourn\'ees
  {S}oc. {M}ath. {F}rance, {I}nst. {R}echerche {M}ath. {A}vanc\'ee,
  {S}trasbourg, 1973)}, pages 35--59. Lecture Notes in Math., Vol. 404.
  Springer, Berlin, 1974.

\bibitem{MR1071805}
Persi Diaconis and James~Allen Fill.
\newblock Strong stationary times via a new form of duality.
\newblock {\em Ann. Probab.}, 18(4):1483--1522, 1990.

\bibitem{MR2530103}
Persi Diaconis and Laurent Miclo.
\newblock On times to quasi-stationarity for birth and death processes.
\newblock {\em J. Theoret. Probab.}, 22(3):558--586, 2009.

\bibitem{MR0193671}
E.~B. Dynkin.
\newblock {\em Markov processes. {V}ols. {I}, {II}}, volume 122 of {\em
  Translated with the authorization and assistance of the author by J. Fabius,
  V. Greenberg, A. Maitra, G. Majone. Die Grundlehren der Mathematischen
  Wissenschaften, B\"ande 121}.
\newblock Academic Press Inc., Publishers, New York, 1965.

\bibitem{MR0328162}
G.~K. Eagleson.
\newblock A characterization theorem for positive definite sequences on the
  {K}rawtchouk polynomials.
\newblock {\em Austral. J. Statist.}, 11:29--38, 1969.

\bibitem{MR2530102}
James~Allen Fill.
\newblock The passage time distribution for a birth-and-death chain: strong
  stationary duality gives a first stochastic proof.
\newblock {\em J. Theoret. Probab.}, 22(3):543--557, 2009.

\bibitem{MR0284628}
George Gasper.
\newblock Positivity and the convolution structure for {J}acobi series.
\newblock {\em Ann. of Math. (2)}, 93:112--118, 1971.

\bibitem{MR0310536}
George Gasper.
\newblock Banach algebras for {J}acobi series and positivity of a kernel.
\newblock {\em Ann. of Math. (2)}, 95:261--280, 1972.

\bibitem{MR2654550}
Laurent Miclo.
\newblock On absorption times and {D}irichlet eigenvalues.
\newblock {\em ESAIM Probab. Stat.}, 14:117--150, 2010.

\bibitem{miclo:hal-01117051}
Laurent Miclo.
\newblock {On the hypergroup property}.
\newblock Preprint available at \texttt{https://hal.archives-}\\
  \texttt{ouvertes.fr/hal-01117051}, February 2015.

\end{thebibliography}
\end{document}